\documentclass{amsart}
\usepackage{amssymb,amsfonts, color,epsf}
\usepackage [cmtip,arrow]{xy}
\xyoption{all}
\usepackage {pb-diagram,pb-xy}
\usepackage{enumerate}
\usepackage{accents}
\usepackage[noautoscale]{youngtab}

\ifx\pdfpageheight\undefined
   \usepackage[dvips,colorlinks=true,linkcolor=blue,citecolor=red,%
      urlcolor=green]{hyperref}
   \usepackage[dvips]{graphicx}
   \makeatletter
   \edef\Gin@extensions{\Gin@extensions,.mps}
   \makeatother
\else
 \usepackage[pdftex]{graphicx}
   \usepackage[bookmarksopen=false,pdftex=true,breaklinks=true,%
      backref=page,pagebackref=true,plainpages=false,%
      hyperindex=true,pdfstartview=FitH,colorlinks=true,%
      pdfpagelabels=true,colorlinks=true,linkcolor=blue,%
      citecolor=red,urlcolor=green,hypertexnames=false%
      ]%
   {hyperref}
\fi
\usepackage{bm}

\usepackage{amsmath,amssymb,amsfonts}
\newtheorem{theorem}{Theorem}
\newtheorem{lemma}{Lemma}

\newtheorem{proposition}{Proposition}
\newtheorem{conjecture}{Conjecture}

\newtheorem{question}{Question}

\theoremstyle{definition}
\newtheorem{definition}{Definition}
\newtheorem{example}{Example}

\newtheorem{notation}{Notation}

\theoremstyle{remark}
\newtheorem{remark}{Remark}

\definecolor{DarkBlue}{rgb}{0,0.1,0.55}

\numberwithin{equation}{section}

\newcommand {\hide}[1]{}
 \newcommand {\sign} {\mbox{\bf sign}}

\newcommand {\junk}[1]{}

\newcommand {\R} {\mathrm{R}}

\newcommand {\C}     {\mathrm{C}}

\newcommand {\Sphere}{\mbox{${\bf S}$}}     

\newcommand {\Z}  {\mathbb{Z}}
 \newcommand {\N}         {\mathbb{N}}
\newcommand {\Q}         {\mathbb{Q}}
\newcommand {\kk}         {\mathbf{k}}
\newcommand{\kbb}    {\mathbb{K}}

\newcommand {\mm}         {\mathbf{m}}
\newcommand{\dd} {\mathbf{d}}

\newcommand{\F}{\mathbb{F}}

\newcommand {\ZZ} {{\rm Z}}
\newcommand {\RR} {{\mathcal R}}

\newcommand {\la}   {{\langle}}
\newcommand {\ra}   {{\rangle}}
\newcommand {\eps} {{\varepsilon}}

\newcommand {\Ker}      {\mbox{\rm ker}}

\newcommand {\PP}     {\mathbb{P}} 

\newcommand{\card}{\mathrm{card}}
\newcommand{\rank}{\mathrm{rank}}

\def\addots{\mathinner{\mkern1mu
\raise1pt\vbox{\kern7pt\hbox{.}}
\mkern2mu\raise4pt\hbox{.}\mkern2mu
\raise7pt\hbox{.}\mkern1mu}}

\newcommand{\HH}  {\mbox{\rm H}}

\newcommand{\Hom}{\textrm{hom}}

\newcommand{\x}{\mathbf{x}}
\newcommand{\X}{\mathbf{X}}

\newcommand{\y}{\mathbf{y}}
\newcommand{\Y}{\mathbf{Y}}

\newcommand{\ZB}{\mathbf{Z}}

\newcommand{\Par}{\mathrm{Par}}

\newcommand{\Ind}{\mathrm{Ind}}

\newcommand{\Def}{\mathrm{Def}}

\newcommand{\Ext}{\mathrm{Ext}}
\newcommand{\length}{\mathrm{length}}

\newcommand{\gdom}{{\;\underline{\triangleright}\;}}

\newcommand{\grad}{\mathrm{grad}}
\newcommand{\ind}{\mathrm{ind}}
\newcommand{\orbit}{\mathrm{orbit}}
\newcommand{\mult}{\mathrm{mult}}

\newcommand{\m}{\mathbf{m}}

\begin{document}
\title[On the isotypic decomposition of cohomology modules]
{
On the isotypic decomposition of cohomology modules of symmetric semi-algebraic sets: polynomial bounds on multiplicities
}
\author{Saugata Basu}
\address{Department of Mathematics,
Purdue University, West Lafayette, IN 47906, U.S.A.}
\email{sbasu@math.purdue.edu}

\author{Cordian Riener}
\address{Department of Mathematics and Statistics, Faculty of Science and Technology, University of Troms\o, 9037 Troms\o,Norway}
\address{Fachbereich Mathematik und Statistik, Universit\"at Konstanz, 78457 Konstanz, Germany}

\email{cordian.riener@uni-konstanz.de}


\date{\textbf{\today}}

\thanks{Part of this research was performed while the authors were visiting the Institute for Pure and Applied Mathematics (IPAM), which is supported by the National Science Foundation.
Basu was also  partially supported by NSF grants
CCF-1319080, CCF 1618981, DMS-1161629,  and DMS-1620271.
 }

\begin{abstract}
We consider symmetric 
(under the action of products of finite symmetric groups)
real algebraic varieties and semi-algebraic sets, as well
as symmetric complex varieties in affine and projective spaces,
defined by polynomials of  degrees bounded by a fixed constant $d$. 
We prove that if a Specht module, $\mathbb{S}^\lambda$, 
appears with positive multiplicity in the isotypic decomposition of the cohomology
modules of such sets, then the rank of the partition $\lambda$ is bounded by $O(d)$. This implies a 
polynomial (in the dimension of the ambient space) bound on the number of such modules.  Furthermore, we prove a polynomial
bound on the multiplicities of  those that do appear  with positive multiplicity
in the isotypic decomposition of the above mentioned cohomology modules.

We  give some applications of our methods in proving lower bounds on the degrees
of defining polynomials of certain symmetric semi-algebraic sets, as well as improved bounds on the
Betti numbers of the images under projections of (not necessarily symmetric) bounded
real algebraic sets, improving in certain situations prior results of Gabrielov, Vorobjov and Zell.
\hide{
We conjecture that the multiplicities 
of the irreducible representations of the symmetric group  
in the cohomology modules of symmetric semi-algebraic sets defined by polynomials having fixed degrees
are computable with polynomial complexity, which would imply
that the Betti numbers of such sets are also computable with polynomial complexity. This is in contrast with general semi-algebraic sets, for which this problem is provably hard 
($\#\mathbf{P}$-hard). 
We also formulate a  representational stability question, 
asking whether the multiplicities associated of the Specht-modules associated to any fixed partition in the cohomology 
sequences of varieties  $(V_k(I))_{k >0}$ defined by a  finitely generated symmetric ideal $I$ of the polynomial ring in infinitely many variables 
is ultimately a polynomial in $k$. A positive answer will lead to interesting invariants of such ideals.
}
 \end{abstract}

\maketitle
\tableofcontents

\section{Introduction}
\label{sec:intro}

For any Lie group $G$,  a real or complex variety $V$ equipped with a $G$-action, and a field of coefficients $\F$,the cohomology groups, $\HH^*(V,\F)$, 
of $V$ inherit a structure of a $G$-module. In this paper,  we consider the special case when $G$ is a finite group, and more specifically a product of symmetric groups, $\mathfrak{S}_\kk = \mathfrak{S}_{k_1} \times \cdots \times \mathfrak{S}_{k_\omega}$, 
acting linearly on finite dimensional real and complex vector spaces by the standard action of permuting coordinates, and $\F$ a field of characteristic $0$.
(Note that
 the topological structure of varieties (also symmetric spaces) 
admitting actions of Lie groups is a very well-studied
topic (see for example \cite{Mimura-Toda-book}). Here we concentrate on the action of finite reflection groups,
which seems to be a less developed field of study.)
We study quantitatively, the $\mathfrak{S}_\kk$-module structure of the cohomology groups of $\mathfrak{S}_\kk$-symmetric algebraic varieties,
and more generally semi-algebraic sets.
We prove upper bounds on the multiplicities of the various irreducibles that appear in the isotypic decomposition
of these modules, as well as restrictions on those that are allowed to appear with non-zero multiplicities. Our upper bounds (both on the multiplicities as well as on the number of irreducibles that are allowed) are polynomial in the number of variables, as long as the degrees of the polynomials defining the variety or semi-algebraic set are held fixed.  We give a couple of applications of these results in proving lower bounds on degrees, as well as improving existing bounds on the Betti numbers of images of semi-algebraic sets (not necessarily symmetric) under polynomial maps. 

We begin  with some history and motivation behind studying these questions.

\subsection{History and motivation}
\label{subsec:history}
Throughout this paper $\R$ will denote a fixed real closed field and $\C$ the algebraic closure of $\R$. 
We also fix a field $\F$ of characteristic $0$.
For any 
closed
semi-algebraic set $S$ we will denote by $b^i(S,\F)$ the dimension of the $i$-th cohomology group,
$\HH^i(S,\F)$, and by $b(S,\F) = \sum_{i \geq 0} b^i(S,\F)$.
 (We refer the reader to 
\cite[Chapter 6]{BPRbook2} for the definition of homology/cohomology  groups of semi-algebraic
sets defined over arbitrary real closed fields, noting that they are isomorphic to the singular 
homology/cohomology groups in the special case of $\R = \mathbb{R}$.)

\subsection{Non-equivariant bounds}
The problem of obtaining quantitative bounds on the topology measured by the the Betti numbers 
of real semi-algebraic as well as
complex constructible sets in terms of the degrees and the number of defining polynomials is
very well studied (see for example, \cite{BPR10} for a survey). For semi-algebraic (respectively, constructible) subsets of $\R^k$ (respectively, $\C^k$) defined by $s$ polynomials of degrees bounded by $d$, these bounds are typically exponential in $k$, and polynomial (for fixed $k$) in $s$ and $d$.

More precisely, suppose that $S$  is a semi-algebraic 
(resp. constructible) subset of $\R^k$  (resp. $\C^k$) defined by a quantifier-free formula involving
$s$ polynomials in $\R[X_1,\ldots,X_k]$ (resp. $\C[X_1,\ldots,X_k]$) of degrees bounded by $d$.

\begin{theorem}[Ole{\u\i}nik and Petrovski{\u\i} \cite{OP}, Thom \cite{T}, Milnor \cite{Milnor2}, \cite{GV07}]
\label{thm:classical}
\[
b(S,\F) \leq (s k d)^{O(k)}.
\]
\end{theorem}

The single exponential dependence on $k$ of the bound in Theorem \ref{thm:classical} is unavoidable. In the real case it suffices to consider the real variety
\begin{equation}
\label{eqn:basic}
V_k = \{1,\ldots,d\}^k \subset \R^k
\end{equation}
defined by the polynomial
\[
F_k = \sum_{i=1}^k \prod_{j=1}^d (X_i-j)^2.
\]
It is easy to see that $\deg(F_k)  =2d$, and  $b_0(V_k) = d^k$.

In the complex case, it follows from a classical formula of algebraic geometry \cite{Hirzebruch-book}
that the sum of the Betti numbers of a non-singular hypersurface  $V_k \subset \PP_\C^k$
of degree $d$ is asymptotically $\Theta(d)^k$. 
Using a standard excision argument and induction on
dimension, the same
asymptotic estimate on the Betti numbers hold for the affine part of such a variety as well.

\subsection{Motivation for studying the equivariant case}
\label{subsec:algorithmic-motivation}
The problem of obtaining tighter estimates on the Betti numbers of semi-algebraic sets 
(motivated partly by applications in other areas of mathematics and theoretical computer science)
has been considered by several authors \cite{Basu1,GV07,BPR8}. 
The algorithmic problem of designing efficient
algorithms for computing these invariants has attracted attention as well \cite{BPRbettione,Bas05-first}. Most of this
work has concentrated on the real semi-algebraic case (since using the real structure a complex
constructible subset $S \subset \C^k$ can be considered as a real semi-algebraic subset of $\R^{2 k}$
defined by twice as many polynomials of the same degrees as those defining $S$), but the complex
case has also being considered separately as well \cite{Scheiblechner07, Walther1}. From the point of
view of algorithmic complexity, the problem of computing the Betti numbers is provably a hard problem -- and so in its full generality a polynomial time algorithm for solving this problem is not to be expected, except in special situations (see \cite{BP'R07joa, Bas05-top} for some of these exceptional
cases). However, an algorithm with even a singly exponential complexity is not known for computing all
the Betti numbers. 

It is a (unproven)  \emph{meta-theorem}  in algorithmic semi-algebraic geometry -- that the worst-case topological complexity of a class of semi-algebraic sets
(measured by the Betti numbers for example)
serve as a rough lower bound for the complexity of algorithms for computing topological
invariants or deciding topological properties of this class of sets.
So the best complexity known for algorithms for determining whether a general semi-algebraic set is empty or connected  is singly exponential, reflecting the singly exponential behavior of the topological
complexity of such sets as exhibited by the example given in Example \ref{eg:basic}. This is true even if the degrees of the
polynomials describing the given set
is bounded by some constant $> 2$). 
On the other hand there are certain classes of semi-algebraic sets where the situation is better.
For example, for semi-algebraic sets defined by few (i.e. any constant number of) quadratic inequalities, we have polynomial upper bounds on the Betti numbers \cite{Bar97}, as well as algorithms 
with polynomial complexities for computing them \cite{BP'R07joa}.

It is intuitively clear that the symmetry imposes strong restrictions on the topology of such sets.
Nevertheless, as shown in Example \ref{eg:basic} below,  the Betti numbers of such sets can be exponentially large.  
However, when the degrees of the defining polynomials are fixed, a polynomial bound is proved on the \emph{equivariant} Betti numbers of such sets in \cite{BC-advances}. 
(These bounds have been subsequently tightened using different methods in \cite{BC-duke}, but these
tighter estimates are not relevant for the current paper.)

On the algorithmic side,
an algorithm with polynomially bounded complexity is
given in \cite{BC-conm} for computing the (generalized) Euler-Poincar\'e characteristics of symmetric 
semi-algebraic sets and their quotients by the action of the symmetric group using techniques developed in \cite{BC-advances}. 
An algorithm with polynomially bounded complexity for computing the Betti numbers of the quotients of such
sets is given in \cite{BC-duke}.
 
Thus, from the point of view of the \emph{meta-theorem} mentioned above, symmetric
semi-algebraic sets pose a dilemma. On one hand their Betti numbers can be exponentially large in the worst case, on the other hand there are reasons to believe that their topological invariants
(when the degree is fixed) has some structure allowing for efficient computation. The polynomial
bound on the equivariant Betti numbers proved in \cite{BC-advances} is the first indication of such a
structure.

\subsection{Summary of the main contributions}
\label{subsec:summary}
We summarize here the main contributions of the current paper. 

\begin{enumerate}[1.]
\item
We consider real as well as complex varieties, and semi-algebraic sets,  on which a product of symmetric
groups acts linearly permuting coordinates. This setting is similar to, but more general than that considered in
\cite{BC-advances,  BC-duke} in that we let the symmetric group act by permuting blocks of variables at a time 
(in \cite{BC-advances, BC-duke} the size of such blocks was limited to one). 
This extra generality is essential in some applications  (see below). The key technical result which makes this generality possible
is Proposition \ref{prop:half-degree}, which generalizes similar results in 
\cite{BC-advances, Riener,Timofte03}  to blocks of sizes larger than one (see also \cite{GRW2016} for an algorithmic application of this result).
\item
Instead of studying the cohomology of the  quotients, $V/\mathfrak{S}_k$, 
where $V$ is a symmetric real or complex variety, or a semi-algebraic set in $k$-dimensional affine or projective space, 
we study the isotypic decomposition of the  $\mathfrak{S}_k$-module  $\HH^*(V,\F)$, where $\F$ is a field of characteristic $0$. The 
Betti numbers of the quotients, i.e. $\dim_\F \HH^*(V/\mathfrak{S}_k,\F)$, can then be recovered from
the multiplicity of the trivial representation in the isotypic decomposition of $\HH^*(V,\F)$ (which is also the dimension of the invariant subspace $\HH^*(V,\F)^{\mathfrak{S}_k}$). We prove (Theorems \ref{thm:main-product-of-symmetric-quantitative}, \ref{thm:main-product-of-symmetric-quantitative-complex}, \ref{thm:main-product-of-symmetric-sa-quantitative})
polynomial bounds on the multiplicities of \emph{all} irreducibles in the isotypic  decomposition, thus generalizing in the results
in \cite{BC-advances,  BC-duke} where polynomial bounds were proved only on the dimension of the trivial representation.
Note that unlike the trivial representation which is of dimension one, the other irreducible
representations of $\mathfrak{S}_\kk$ can have dimensions which are exponentially large (as is unavoidable since the 
dimension of $\HH^*(V,\F)$ can be exponentially large as in Example \ref{eg:basic}).

Moreover, we prove (see Remark \ref{rem:restriction}) that
the number of irreducibles that are allowed to appear is polynomially bounded (and hence a negligible
fraction as $k \rightarrow \infty$) of all irreducibles (which are in bijection with the set of partitions of $k$).
Thus, 
while the Betti numbers of symmetric semi-algebraic sets can be exponentially large, they can be expressed as a sum of polynomially many numbers (the dimensions of the isotypic components),
and each of these numbers is a product of a multiplicity (which is polynomially bounded) and the
dimension of a Specht module (which can be exponentially large, but efficiently computable due to the hook formula (cf. Theorem \ref{thm:hook})). 

\item
In the special case of the multiplicity of the trivial representations, or equivalently the Betti numbers of the quotients, the bounds proved in the current paper still generalizes those in \cite{BC-advances}, since we consider more general actions (permuting blocks of size greater than one). This extra generality is useful in several applications, and we give two applications. In the first application, this added flexibility allows us to treat the case of symmetric complex projective varieties (Theorem \ref{thm:main-product-of-symmetric-quantitative-complex}), with the
symmetric group permuting blocks of size $2$ (the real and imaginary parts). 
Secondly, we are able to generalize a result in \cite{BC-advances} on bounding the Betti numbers of the image
under projection of a real variety, from the case considered in \cite{BC-advances} where the projection
was along one variable, to more general projections (Theorem \ref{thm:descent2-quantitative-new}). 
The crucial new  ingredient is the generalization of the results 
in \cite{BC-advances} to the case of block size greater than one.
\item
Finally, we ask a question and make a conjecture suggested by the results in this paper. 
The question
(Question \ref{question:stability})  is motivated by similar representational stability results in the theory of finitely generated FI-modules \cite{Church-et-al} and asks whether the multiplicities of the irreducible corresponding to some fixed partition should ultimately stabilize to a polynomial for certain naturally defined sequences of varieties. 
We also make the algorithmic conjecture (Conjecture \ref{conj:poly}), stating that the ordinary Betti numbers of symmetric varieties defined by polynomials of fixed degrees should be computable with polynomially bounded complexity. 
We give some evidence in favor of these conjectures.   
\end{enumerate}

\begin{remark}[Homology versus cohomology]
\label{rem:homology}
Note that since $\F$ is a field of characteristic $0$, $\HH^*(S,\F) \cong \Hom(\HH_*(S,\F),\F)$ as vector spaces.
Moreover, from the basic property of $\mathfrak{S}_k$ that the conjugacy class of an element equals that of its inverse it follows that for any finite-dimensional representation $W$ of $\mathfrak{S}_k$, $\Hom(W,\F) \cong W$ as $\mathfrak{S}_k$-modules.
Taken together this implies that 
$\HH_*(S,\F),\HH^*(S,\F)$ for any symmetric semi-algebraic set $S \subset \R^k$, are isomorphic
as  $\mathfrak{S}_k$-modules, and thus for the purposes of  determining the multiplicities of irreducible representations it does not matter whether we consider homology or cohomology modules.
\end{remark}

\hide{
The equivariant
cohomology group, $\HH^*_{\mathfrak{S}_k}(S,\F)$, is isomorphic to the subspace of
$\HH^*(S,\F)^{\mathfrak{S}_k}$ that is fixed by $\mathfrak{S}_k$. Another description of   $\HH^*(S,\F)^{\mathfrak{S}_k}$ is that it  is the isotypic component of $\HH^*
(S,\F)$ corresponding to the trivial one-dimensional representation of $\mathfrak{S}_k$, and the
dimension of $\HH^*(S,\F)^{\mathfrak{S}_k}$ is equal to the multiplicity of the trivial representation
in $\HH^*(S,\F)$. Thus, the main result of \cite{BC-advances} can be expressed as saying that the multiplicity of the trivial representation in   $\HH^*(S,\F)$ is bounded polynomially. It is well known
(see \S \ref{subsec:representation-of-Sn} below) that the irreducible representations of $\mathfrak{S}_k$ (the so called 
\emph{Specht-modules}) are in correspondence with the finite set of \emph{partitions} of $k$. The number of partitions of $k$ is exponentially large. 
}
\subsection{Basic notation and definition}
\label{subsec:basic-notation-definition}
In this section we introduce notation and definitions that we will use for the rest of the paper.

\begin{notation}[Zeros]
\label{not:zeros}
  For $P \in \R [X_{1} , \ldots ,X_{k}]$ (respectively $P \in \C [ X_{1} ,
  \ldots ,X_{k} ]$) we denote by $\ZZ(P, \R^{k})$ (respectively
  $\ZZ (P, \C^{k})$) the set of zeros of $P$ in
  $\R^{k}$(respectively $\C^{k}$). More generally, for any finite set
  $\mathcal{P} \subset \R [ X_{1} , \ldots ,X_{k} ]$ (respectively
  $\mathcal{P} \subset \C [ X_{1} , \ldots ,X_{k} ]$), we denote by $\ZZ
  (\mathcal{P}, \R^{k})$ (respectively $\ZZ (\mathcal{P},
  \C^{k})$) the set of common zeros of $\mathcal{P}$ in
  $\R^{k}$(respectively $\C^{k}$).  For a homogeneous polynomial $P \in \R [X_{0} , \ldots ,X_{k-1}]$ (respectively $P \in \C [ X_{0} ,
  \ldots ,X_{k-1} ]$)   we denote by $\ZZ(P, \PP_\R^{k-1})$ (respectively
  $\ZZ (P, \PP_\C^{k-1})$) the set of zeros of $P$ in
  $\PP_\R^{k-1}$(respectively $\PP_\C^{k-1}$). And, more generally, for any finite set of 
  homogeneous polynomials
  $\mathcal{P} \subset \R [ X_{0} , \ldots ,X_{k-1} ]$ (respectively
  $\mathcal{P} \subset \C [ X_{0} , \ldots ,X_{k-1} ]$), we denote by $\ZZ
  (\mathcal{P}, \PP_\R^{k-1})$ (respectively $\ZZ (\mathcal{P},
  \PP_\C^{k-1})$) the set of common zeros of $\mathcal{P}$ in
  $\PP_\R^{k-1}$ (respectively $\PP_\C^{k-1}$). 
\end{notation}

\begin{notation}[Sign conditions, realizations, $\mathcal{P}$- and $\mathcal{P}$-closed semi-algebraic sets]
  \label{not:sign-condition} For any finite family of polynomials $\mathcal{P}
  \subset \R [ X_{1} , \ldots ,X_{k} ]$, we call an element $\sigma \in \{
  0,1,-1 \}^{\mathcal{P}}$, a \emph{sign condition} on $\mathcal{P}$. For
  any semi-algebraic set $Z \subset \R^{k}$, and a sign condition $\sigma \in
  \{ 0,1,-1 \}^{\mathcal{P}}$, we denote by $\RR (\sigma ,Z)$ the
  semi-algebraic set defined by $$\{ \mathbf{x} \in Z \mid \sign (P (
  \mathbf{x})) = \sigma (P)  ,P \in \mathcal{P} \},$$ and call it the
  \emph{realization} of $\sigma$ on $Z$. More generally, we call any
  Boolean formula $\Phi$ with atoms, $P \{ =,>,< \} 0, P \in \mathcal{P}$, to
  be a \emph{$\mathcal{P}$-formula}. We call the realization of $\Phi$,
  namely the semi-algebraic set
  \begin{eqnarray*}
    \RR \left(\Phi , \R^{k} \right) & = & \left\{ \mathbf{x} \in \R^{k} \mid
    \Phi (\mathbf{x}) \right\}
  \end{eqnarray*}
  a \emph{$\mathcal{P}$-semi-algebraic set}. Finally, we call a Boolean
  formula without negations, and with atoms $P \{\geq, \leq \} 0$, $P\in \mathcal{P}$, to be a 
  \emph{$\mathcal{P}$-closed formula}, and we call
  the realization, $\RR \left(\Phi , \R^{k} \right)$, a \emph{$\mathcal{P}$-closed
  semi-algebraic set}.
\end{notation}

The notion of \emph{partitions} of a given integer will play an important role in what follows,
which necessitates the following notation that we fix for the remainder of the paper.

\begin{notation}[Partitions]
  \label{not:Partition1} 
  We denote by $\Par(k)$ the set of \emph{partitions} of $k$, where each partition $\pi \in \Par(k)$
  (also denoted $\lambda \vdash k$) 
  is a tuple  $(\pi_{1} , \pi_{2} , \ldots
  , \pi_{\ell})$, with $\pi_{1} \geq \pi_{2} \geq \cdots \geq
  \pi_{\ell} \geq 1$, and $\pi_{1} + \pi_{2} + \cdots + \pi_{\ell} =k$. We
  call $\ell$ the length of the partition $\pi$, and denote
  $\length(\pi) = \ell$. 
    
  More generally, for any tuple $\mathbf{k}= (k_{1} , \ldots ,k_{\ell})
  \in \Z_{>0}^\ell$, we will denote by
  $\Par(\mathbf{k}) = \Par(k_{1}) \times \cdots
  \times \Par(k_{\ell})$, and for each $\pmb{\pi}= (
  \pi^{(1)} , \ldots , \pi^{(\ell)}) \in
  \Par(\mathbf{k})$, we denote by
  $\length(\pmb{\pi}) = \sum_{i=1}^{\ell} \length(\pi^{(i)})$. 
  We also denote for each $\mathbf{p}= (p_{1} , \ldots , p_{\ell}) \in
  \N^{\ell}$,
  \begin{eqnarray*}
    | \mathbf{p} | & = & p_{1} + \cdots + p_{\ell},\\
    F (\mathbf{k},\mathbf{p}) & = & \card (\{
    \pmb{\pi}= (\pi^{(1)} , \ldots , \pi^{(\ell)})
    \mid \length (\pi^{(i)}) = p_{i} ,  1
    \leq i \leq \ell \}) .
  \end{eqnarray*}
\end{notation}

\begin{notation}[Transpose of a partition and partitions of bounded lengths]
\label{not:Partition2}
For a partition $\lambda =(\lambda_1,\ldots,\lambda_\ell) \vdash k$, we will denote by $\tilde{\lambda}$ the \emph{transpose} of $\lambda$.
More precisely,
$\tilde{\lambda} = (\tilde{\lambda}_1, \ldots,\tilde{\lambda}_{\tilde{\ell}})$, where $\tilde{\lambda}_j = \card(\{i \mid \lambda_i \geq j \})$.
For $k,d\geq 0$, we denote 
\[
\Par(k,d) := \{ \lambda \in \Par(k)  \mid \length(\lambda) \leq d \}.
\]
More generally,  for $\kk = (k_1,\ldots,k_\ell), \dd = (d_1,\ldots,d_\ell)$ we denote 
\[
\Par(\kk,\dd) := \{ \pmb{\lambda} = (\lambda^{(1)},\ldots,\lambda^{(\ell)}) \mid \lambda^{(i)} \in \Par(k_i),  \length(\lambda^{(i)}) \leq d_i, 1 \leq i \leq \ell \}.
\]
When $\dd = (d,\ldots,d)$, we will also use $\Par(\kk,d)$ to denote $\Par(\kk,\dd)$.
\end{notation}

\begin{notation}[Products of symmetric groups]
\label{not:symmetric-group}
For each $k \in \N$, we denote by $\mathfrak{S}_k$ the symmetric group on $k$ letters (or equivalently the Coxeter group $A_{k-1}$).
For $\kk=(k_1,\ldots,k_\ell) \in \Z_{>0}^\ell$ we denote by $\mathfrak{S}_\kk$ the product group 
$\mathfrak{S}_{k_1} \times \cdots \times \mathfrak{S}_{k_\ell}$, and we will
usually denote $k = |\kk| = \sum_{i=1}^{\ell} k_i$.
\end{notation}

We first need some more notation.

\begin{notation}[Young subgroups of product of symmetric groups]
For 
\[\lambda =(\lambda_1,\ldots, \lambda_d) \in \Par(k),
\] 
we will denote by 
$\mathfrak{S}_{\pmb{\lambda}} \cong \mathfrak{S}_{\lambda_1} \times \cdots \times \mathfrak{S}_{\lambda_d}$ the subgroup of 
$\mathfrak{S}_k$ which is the direct product of the subgroups $G_i \cong \mathfrak{S}_{\lambda_i}$, where $G_i$ is the subgroup
of permutations of $[1,k]$ fixing $[1,k] \setminus [\lambda_1+\cdots+\lambda_{i-1}+1, \lambda_1+\cdots+\lambda_i]$.

More generally, for $\kk = (k_1,\ldots,k_\ell) \in \Z_{>0}^\ell$, $\pmb{\lambda} = (\lambda^{(1)},\ldots,\lambda^{(\ell)}) \in \Par(\kk)$, we denote by 
$\mathfrak{S}_{\pmb{\lambda}}$ the subgroup $\mathfrak{S}_{\lambda^{(1)}} \times \cdots \times \mathfrak{S}_{\lambda^{(\ell)}}$ of
$\mathfrak{S}_\kk$, where for $1 \leq j \leq \ell$, $\mathfrak{S}_{\lambda^{(j)}}$ is the subgroup of $\mathfrak{S}_{k_j}$ defined
above. 
\end{notation}

\begin{notation}[Irreducible representations (Specht modules)  of symmetric groups]
For $\lambda \in \Par(k)$, we will denote by $\mathbb{S}^{\lambda}$ the irreducible representation  (over the field $\F$) of 
$\mathfrak{S}_k$ corresponding to
$\lambda$ (see \cite{Procesi-book} for definition). 
Note that $\mathbb{S}^{(k)}$ is the trivial representation (corresponding to the  partition $(k) \in \Par(k)$
(which we also denote by $\mathbf{1}_{\mathfrak{S}_k}$), and 
$\mathbb{S}^{(1^k)}$ is the sign representation, which we will also denote by $\mathbf{sign}_k$.
It is  well known fact that for any $\lambda \in \Par(k)$, 
\[
\mathbb{S}^{(\tilde{\lambda})} \cong \mathbb{S}^{(\lambda)} \otimes \mathbf{sign}_k.
\]

For $\kk = (k_1,\ldots,k_\ell) \in \Z_{>0}^\ell$, $\pmb{\lambda} = (\lambda^{(1)},\ldots,\lambda^{(\ell)}) \in \Par(\kk)$, we denote by $\mathbb{S}^{\pmb{\lambda}}$ the irreducible representation 
$\mathbb{S}^{\lambda^{(1)}} \boxtimes \cdots \boxtimes \mathbb{S}^{\lambda^{(\ell)}}$ of $\mathfrak{S}_{\kk}$.
\end{notation}

\begin{definition}[$\mathfrak{S}_{\kk}$-symmetric polynomials]
\label{def:action}
Let $\kbb$ be the field $\R$ or $\C$.
Suppose that $\kk=(k_1,\ldots,k_\ell), \m =(m_1,\ldots,m_\ell) \in \Z_{>0}^\ell$, and let 
$P \in \kbb[\X^{(1)},\ldots,\X^{(\ell)}]$
where for $1 \leq h \leq \ell$, 
$\X^{(h)} = \left(X^{(h)}_{i,j}\right)_{1\leq i \leq k_h, 1\leq j \leq m_h}$.

The group $\mathfrak{S}_{\kk}$, acts on $\kbb[\X^{(1)},\ldots,\X^{(\ell)}]$
by permuting for each $i,1\leq i \leq \ell$, 
the rows of $\X^{(h)}$ by the group $\mathfrak{S}_{k_h}$. For $\pmb{\pi} \in \mathfrak{S}_\kk$,
and $P \in \kbb[\X^{(1)},\ldots,\X^{(\ell)}]$, we denote the by $\pmb{\pi}\cdot P$ the image of
$P$ under $\pmb{\pi}$. 
We say that \emph{$P$ is $\mathfrak{S}_\kk$-symmetric} if it is invariant under the action of 
$\mathfrak{S}_\kk$, i.e. if $\pmb{\pi}\cdot P = P$ for every $\pmb{\pi} \in \mathfrak{S}_\kk$.

For $\dd =(d_1,\ldots,d_\ell) \in \Z_{> 0}^\ell$, we will denote by $\kbb[\X^{(1)},\ldots,\X^{(\ell)}]^{\mathfrak{S}_\kk}_{\leq \dd}$, the finite dimensional subspace of 
$\kbb[\X^{(1)},\ldots,\X^{(\ell)}]$ consisting of $\mathfrak{S}_\kk$-symmetric 
polynomials whose degree in $\X^{(i)}$ is bounded by $d_i$ for $1 \leq i \leq \ell$.

Similarly, we say that a 
subset $S \subset \kbb^K, K = \sum_{1\leq i \leq \ell} k_i m_i$,
is $\mathfrak{S}_\kk$-symmetric if it is stable under the above action of $\mathfrak{S}_\kk$. 

When $\ell=1,m_1=1$, and $K= k_1m_1 = k$, the action defined above is the usual action of $\mathfrak{S}_k$ on $\kbb^k$
permuting coordinates. 
\end{definition}

\begin{remark}
\label{rem:complex-action}
Note in case $\kbb = \C$, the action of $\mathfrak{S}_\kk$ on $\C^K$ defined above in Definition
\ref{def:action}  can also be seen as the action of $\mathfrak{S}_\kk$ on $\R^{2 K}$ (considering $\C = \R \oplus  i\R$), replacing $\m$ by $2\m$.
\end{remark}

\subsection{Basic example}
\label{subsec:basic-example}
Before proceeding further, we discuss an example which is our guiding example for the rest of the paper. While explaining the example we will 
assume a certain familiarity with the representation theory of symmetric groups.
For the convenience of the reader we have included all the facts from the representation theory of symmetric groups that we need 
in \S \ref{subsec:representation-of-Sn} (and which the reader can consult if needed).

\begin{example}[Real affine case]
\label{eg:basic}
Let 
\[
F_k = \sum_{i=1}^k X_i^2(X_i -1)^2 - \eps,
\]
and 
\begin{equation}
\label{eqn:eg:basic}
V_k = \ZZ(F_k,\R^k).
\end{equation}

Then, for all  $\eps, 0 < \eps \ll 1$,
$V_k$ is a closed and bounded non-singular hypersurface in $\R^k$, (in fact also in $\PP_\R^k$), 
the semi-algebraic set $S_k$ defined by $F_k \leq 0$ is homotopy equivalent to the finite set 
of points $\{0,1\}^k$, and is bounded by $V_k$.

Clearly, $b_0(V_k,\F) = 2^k$, and 
it follows from Poincar\'e duality applied to $V_k$  that $b_{k-1}(V_k,\F) = 2^k$ as well.
It also follows from Alexander-Lefshetz duality that $\HH^i(V_k,\F)=0$ for $0< i < k-1$.
 
The real algebraic variety $V_k$ is symmetric under the standard action of the symmetric group $\mathfrak{S}_k$ on $\R^k$ 
permuting the coordinates. This action induces an $\mathfrak{S}_k$-module structure on  $\HH^*(V_k,\F)$,
and it is interesting to study the isotypic decomposition of this representation into its isotypic components corresponding to the
various irreducible representations of $\mathfrak{S}_k$, namely the Specht modules $\mathbb{S}^\lambda$ indexed by
different partitions $\lambda\vdash k$ (see for example \cite{Procesi-book} for the definition of Specht modules).

We now describe this decomposition.

\[
\HH^0(V_k,\F)  \cong \bigoplus_{0 \leq i \leq k} \HH^0(V_{k,i},\F), 
\] 
where for $0 \leq i \leq k$, $V_{k,i}$ is the $\mathfrak{S}_k$-orbit of the 
connected component of $V_k$ infinitesimally close (as a function of $\eps$)  to the
point $\x^i = (\underbrace{0,\ldots,0}_i,\underbrace{1,\ldots,1}_{k-i} )$,
and $\HH^0(V_{k,i},\F)$ is a sub-representation of $\HH^0(V_k,\F)$.
  
It is also clear that the isotropy subgroup of 
the class in $\HH^0(V_k,\F)$ corresponding to $V_{k,i}$
is isomorphic to $\mathfrak{S}_i \times \mathfrak{S}_{k-i}$, and
hence,
\begin{eqnarray*}
\HH^0(V_{k,i},\F) &\cong& \Ind_{\mathfrak{S}_i \times \mathfrak{S}_{k-i}}^{\mathfrak{S}_k}( \mathbb{S}^{(i)} \boxtimes \mathbb{S}^{(k-i)}) \\
&\cong& M^{(i,k-i)} \mbox{ if }  i \geq k-i,\\
&\cong & M^{(k-i,i)} \mbox{ otherwise}.  
\end{eqnarray*}

where for any $\lambda \vdash k$, we denote by $M^\lambda$ the Young module corresponding to
$\lambda$ (see Definition \ref{def:Young}).

Also, observe that $\HH^0(V_{k,i},\F)$ and $\HH^0(V_{k,k-i},\F)$ are isomorphic  as $\mathfrak{S}_k$-modules.  
In the following, for partitions $\mu,\lambda \vdash k$, we will denote by $K(\mu,\lambda)$ the corresponding
\emph{Kostka number} (see Definition \ref{def:Kostka} below). For this example, it is sufficient to observe that
if $\mu \gdom \lambda$ (see Definition \ref{def:dominance} for the definition of the \emph{dominance order} $\gdom$ on the set of
partitions), and if $\mu$ has at most $2$ rows, then $K(\mu,\lambda) = 1$.
It now follows from 
Proposition \ref{prop:Young}
that for $k$ odd,
\begin{eqnarray*}
\HH^0(V_k,\F) &\cong & \bigoplus_{\substack{\lambda \vdash k\\ \ell(\lambda) \leq 2}}  (M^\lambda \oplus M^\lambda) \\
& \cong & \bigoplus_{\substack{\lambda \vdash k\\ \ell(\lambda) \leq 2}} \bigoplus_{\mu \gdom \lambda} 2 K(\mu,\lambda)  \mathbb{S}^\mu\\
& \cong & \bigoplus_{\substack{\lambda \vdash k\\ \ell(\lambda) \leq 2}} \bigoplus_{\mu \gdom \lambda} 2 \mathbb{S}^\mu\\
& \cong & \bigoplus_{\substack{\mu \vdash k\\ \ell(\mu) \leq 2}} m_\mu \mathbb{S}^\mu,
\end{eqnarray*}
where  for each $\mu = (\mu_1,\mu_2)  \vdash k$,   
\begin{eqnarray*}
m_\mu &=& 2(\mu_1 - \lfloor k/2 \rfloor)  \\
            &=&  2\mu_1 - k +1. 
\end{eqnarray*}

For $k$ even we have,
\begin{eqnarray*}
\HH^0(V_k,\F) &\cong & \left(\bigoplus_{\substack{\lambda \vdash k\\ \ell(\lambda) \leq 2 \\ \lambda \neq (k/2,k/2)}}  (M^\lambda \oplus M^\lambda) \right)  \bigoplus M^{(k/2,k/2)}  \\
& \cong & \left(
\bigoplus_{\substack{\lambda \vdash k \\ \ell(\lambda) \leq 2 \\ \lambda \neq (k/2,k/2)}} \bigoplus_{\mu \gdom \lambda} 2 K(\mu,\lambda)  \mathbb{S}^\mu 
\right) \oplus
 \left(
 \bigoplus_{\mu \gdom (k/2,k/2)} K(\mu,(k/2,k/2))\mathbb{S}^\mu
 \right)\\
& \cong & \bigoplus_{\substack{\mu \vdash k \\ \ell(\mu) \leq 2}} m_\mu \mathbb{S}^\mu,
\end{eqnarray*}
where for each $\mu = (\mu_1,\mu_2)  \vdash k$,  
\begin{eqnarray*}
m_\mu &=& 2(\mu_1 - {k/2})+1 \\
	     &=& 2\mu_1 -k +1. 
\end{eqnarray*}

We deduce for all $k$, 
\begin{eqnarray*}
m_\mu &=& 2\mu_1 -k +1 \label{eqn:even-and-odd} \\
            &\leq & k+1. \nonumber
\end{eqnarray*}

For $\mu = (\mu_1,\mu_2) \vdash k$, by the hook-length formula (Eqn. \eqref{eqn:hook}) 
we have,
\begin{eqnarray}
\label{eqn:hook-length}
\dim \; \mathbb{S}^\mu &=& \frac{k! \; (\mu_1 - \mu_2+1)}{(\mu_1+1)!\mu_2!}.
\end{eqnarray}

This completes the description of the isotypic decomposition of $\HH^0(V_k,\F)$. 

In particular for $k=2,3$ we have:
\begin{eqnarray*}
\HH^0(V_2,\F) &\cong &  3\mathbb{S}^{(2)} \oplus \mathbb{S}^{(1,1)}, \\
\HH^0(V_3,\F) &\cong &  4\mathbb{S}^{(3)} \oplus 2\mathbb{S}^{(2,1)}.
\end{eqnarray*}

The isotypic decomposition of $\HH^{k-1}(V_k,\F)$ requires one further ingredient -- namely, an
$\mathfrak{S}_k$-equivariant version of the classical Poincar\'e duality theorem for oriented manifolds. We include 
a proof of this result (Theorem \ref{thm:poincare-duality}) in  \S \ref{subsec:Poincare-duality}.

We note that $V_k$ is a closed and bounded real orientable manifold, 
by Poincar\'e duality theorem there exists an isomorphism between $\HH^0(V,\F)$ and $\HH^{k-1}(V,\F)$. This
isomorphism is not necessarily a $\mathfrak{S}_k$-module isomorphism. However, it follows from  
Theorem \ref{thm:poincare-duality} (which is a stronger form of Poincar\'e duality for orientable symmetric manifolds) 
that the isotypic representation
of $\HH^{k-1}(V_k,\F)$ is isomorphic (as an $\mathfrak{S}_k$-module) to 
$\HH^0(V_k,\F) \otimes
\textbf{sign}_{k}
$.

Thus, denoting for each $\lambda \vdash k$, the \emph{transpose} of the partition $\lambda$ by $\tilde{\lambda}$, 

\begin{eqnarray*}
\HH^{k-1}(V_k,\F) &\cong &  \bigoplus_{\substack{\mu \vdash k\\ \ell(\mu) \leq 2}} m_\mu \mathbb{S}^{\tilde{\mu}},
\end{eqnarray*}
where  for  each $\mu = (\mu_1,\mu_2)  \vdash k$,  
$m_\mu$ is defined above in \eqref{eqn:even-and-odd}.
In particular for $k=2,3$ we have:
\begin{eqnarray*}
\HH^1(V_2,\F) &\cong &  3\mathbb{S}^{(1,1)} \oplus \mathbb{S}^{(2)}, \\
\HH^2(V_3,\F) &\cong &  4\mathbb{S}^{(1,1,1)} \oplus 2\mathbb{S}^{(2,1)}.
\end{eqnarray*}

Notice that the multiplicity $m_{1^k}$ of the Specht module 
$\mathbb{S}^{1^k} = \textbf{sign}_{k}$ in
$\HH^0(V_k,\F)$ is equal to $0$ for $k>2$. This implies that the multiplicity of the trivial
representation $\mathbb{S}^{(k)}$ is equal to $0$ in $\HH^{k-1}(V_k,\F)$, and thus
$\HH^{k-1}_{\mathfrak{S}_k}(V_k,\F)=0$ as well (for $k >2$).

 Also, notice that  the multiplicity of each Specht-module, $\mathbb{S}^\mu, \mu \vdash k$, 
in the isotypic decomposition of $\HH^*(V_k,\F)$ is bounded
polynomially (in fact, linearly) in $k$, but the dimension of $\HH^*(V_k,\F)$ itself  is exponentially large in $k$.  

Note that  since $\dim \HH^0(V_k,\F) =  2^k$, we obtain as a consequence  (from \eqref{eqn:even-and-odd} and 
 \eqref{eqn:hook-length}) the identity
 
 \begin{eqnarray*}
  k!\;\left(\sum_{\substack{
  \mu_1 \geq \mu_2\geq 0\\ \mu_1+\mu_2 =k} }  \frac{(\mu_1 - \mu_2 +1)^2}{(\mu_1+1)!\mu_2!}\right)
   &=& 2^k
  \end{eqnarray*}
  (which can also be proved easily by more elementary means).
\end{example}

\begin{example}[Projective case]
\label{eg:real-projective}
Let 
\[
P = \sum_{0\leq i< j \leq k-1} (X_i^2 - X_j^2)^2, 
\]
and let $W_k = \ZZ(P,\PP_\R^{k-1})$.
Then,
\[
W_k = \{(x_0:\cdots:x_{k-1}) \mid x_i = \pm 1, 0 \leq i \leq k-1\},
\]
and is symmetric under the action of $\mathfrak{S}_{k}$ on $\PP_\R^{k-1}$ permuting the homogeneous coordinates.

It is clear that 
\[
\HH^0(W_k,\F) \cong\HH^0(V_k,\F),
\]
where $V_k$ is the real affine variety defined in \eqref{eqn:eg:basic}, and the stated isomorphism is
an isomorphism of $\mathfrak{S}_k$-modules. 
\end{example}

\subsection{Equivariant cohomology}
\label{subsec:equivariant}
We recall also the definition of \emph{equivariant cohomology groups} of a
$G$-space for an arbitrary compact Lie group $G$. For $G$ any compact Lie
group, there exists a \emph{universal principal $G$-space}, denoted $E G$,
which is contractible, and on which the group $G$ acts freely on the right. 
The \emph{classifying space} $B G$, is the orbit space of this action, i.e.
$B G= E G/G$.

\begin{definition}[Equivariant cohomology]
  \label{def:equivariant-cohomology} (Borel construction) 
  Let $X$ be a space with a left action of  the group $G$. 
  Then, $G$ acts diagonally on the space $E G \times X$ by $g (z,x) = (z \cdot g^{-1} ,g \cdot x)$. For any
  field of coefficients $\F$, the \emph{$G$-equivariant cohomology
  groups of $X$} with coefficients in $\F$, denoted by
  $\HH^{\ast}_{G} (X,\F)$, is defined by
$\HH^{\ast}_{G} (X,\F)  =  \HH^{\ast} (E G \times X/G,\F)$.
\end{definition}

In the situation of interest in the current paper, where $G= \mathfrak{S}_\kk$ acting on
a $\mathfrak{S}_\kk$-symmetric semi-algebraic subset $X \subset \R^k$,
and $\F$ is a field with characteristic equal to $0$,
we have the isomorphisms (see \cite{BC2013}):

\begin{equation}
\label{eqn:iso}
\HH^{\ast} (S/\mathfrak{S}_\kk,\F) \xrightarrow{\sim} \HH_{\mathfrak{S}_\kk}^{\ast} (S,\F
) \xrightarrow{\sim} \HH^\ast(S,\F)^{\mathfrak{S}_\kk}.
\end{equation}

\subsection{Prior work}
\label{subsec:prior-work}
The problem of bounding the equivariant Betti numbers of symmetric semi-algebraic subsets of 
$\R^k$ was investigated in \cite{BC-advances}. We recall in this section a few results from
\cite{BC-advances} that are generalized in the current paper.

We recall some definitions and notation from \cite{BC-advances}.

\begin{notation}[Equivariant Betti numbers]
  \label{not:equivariant-betti} 
  For any $\mathfrak{S}_{\mathbf{k}}$ symmetric
  semi-algebraic subset $S \subset \R^{k}$ with $\mathbf{k}= (k_{1} , \ldots
  ,k_{\ell}) \in \N^{\ell}$, with $k= \sum_{i=1}^{\ell} k_{i}$,
  and any field $\F$,  we denote
  \begin{eqnarray*}
    b_{\mathfrak{S}_{\mathbf{k}}}^{i} (S,\F) & =  & b_{i} (S/\mathfrak{S}_{\mathbf{k}} ,\F) ,\\
    b_{\mathfrak{S}_{\mathbf{k}}} (S,\F) & = & \sum_{i \geq 0}
    b_{\mathfrak{S}_{\mathbf{k}}}^{i} (S,\F).
  \end{eqnarray*}
\end{notation}

The following theorem is proved in \cite{BC-advances}.

\begin{theorem} \cite[Theorem 6]{BC-advances}
  \label{thm:main} Let $\mathbf{k}= (k_{1} , \ldots ,k_{\ell}) \in
  \N^{\ell}$,with  $k= \sum_{i=1}^{\ell} k_{i}$. Suppose that 
  $P \in \R [\X^{(1)} , \ldots,\X^{(\ell)} ]$, where each
  $\X^{(i)}$ is a block of $k_{i}$ variables, is a 
  non-negative polynomial, such that $V= \ZZ(P, \R^{k})$ is
  invariant under the action of $\mathfrak{S}_{\mathbf{k}}$ permuting each
  block $\X^{(i)}$ of $k_{i}$ coordinates. Let
  $\deg_{\X^{(i)}} (P)   \leq  d$ for $1 \leq i \leq \ell$. Then, for any
  field of coefficients $\F$,
  \begin{eqnarray*}
    b (V/\mathfrak{S}_{\mathbf{k}} ,\F) & \leq &
    \sum_{\mathbf{p}= (p_{1} , \ldots , p_{\ell}) ,1 \leq p_{i}
    \leq \min (2d,k_{i})  } F (\mathbf{k},\mathbf{p})  d (2d-1)^{|
    \mathbf{p} | +1}
  \end{eqnarray*}
  (where $F(\mathbf{k},\mathbf{p})$ is defined in Notation \ref{not:Partition1}).
  If for each $i,1 \leq i \leq \ell$, $2d  \leq k_{i}$, then
  \begin{eqnarray*}
    b (V/\mathfrak{S}_{\mathbf{k}} ,\F) & \leq & (k_{1}
    \cdots k_{\ell})^{2d} (O (d))^{2 \ell d+1} .
  \end{eqnarray*}
\end{theorem}

More generally, the following bound holds for symmetric semi-algebraic sets.

\begin{theorem}\cite[Theorem 7]{BC-advances}
  \label{thm:main-sa-closed}
  Let $\mathbf{k}= (k_{1} , \ldots ,k_{\ell})
  \in \N^{\ell}$, with $k= \sum_{i=1}^{\ell} k_{i}$, and let
  $\mathcal{P} \subset \R [ \X^{(1)} , \ldots,\X^{(\ell)} ]$ be a finite set of polynomials,
  where each $\X^{(i)}$ is a block of $k^{(i)}$
  variables, and such that each $P \in \mathcal{P}$ is symmetric in each block
  of variables $\X^{(i)}$. Let $S \subset \R^{k}$
  be a $\mathcal{P}$-closed-semi-algebraic set. Suppose that $\deg (P)  
  \leq  d$ for each $P  \in  \mathcal{P}$, $\card (
  \mathcal{P}) =s$, and let $D=D (\mathbf{k},d) = \sum_{i=1}^{\ell} \min
  (k_{i} ,5d)$. Then, for any field of coefficients $\F$,
  \begin{eqnarray*}
    b (S/\mathfrak{S}_{\mathbf{k}} ,\F) & \leq & \sum_{i=0}^{D-1}
    \sum_{j=1}^{D-i} 
    \binom{2 s+1}{j} 6^{j} G (\mathbf{k},2d)
  \end{eqnarray*}
  where 
    \begin{eqnarray*}
    G (\mathbf{k},d) & = & \sum_{\mathbf{p}= (p_{1} , \ldots ,
    p_{\ell}) ,1 \leq p_{i} \leq \min (2d,k_{i})} F (
    \mathbf{k},\mathbf{p})  d (2d-1)^{| \mathbf{p} | +1}
  \end{eqnarray*}
  (and $F(\mathbf{k},\mathbf{p})$ is defined in Notation \ref{not:Partition1}).
\end{theorem}

\begin{remark}
\label{rem:main-sa-closed}
  In the particular case, when $\ell =1$, $d=O (1)$, the bound in Theorem
  \ref{thm:main-sa-closed} takes the following asymptotic (for $k \gg 1$)
  form.
  \begin{eqnarray*}
    b (S/\mathfrak{S}_{k} ,\F) & \leq & O  (s^{5d} k^{4d-1}) .
  \end{eqnarray*}
\end{remark}

The rest of the paper is organized as follows. In \S \ref{sec:results} we state the new results proved in this paper.
In \S \ref{sec:preliminaries} we prove or recall certain preliminary facts that will be needed in the proofs of the
main theorems.  In \S \ref{sec:proofs-of-main} we prove the main theorems, and finally in \S
\ref{sec:conclusion} we end with some open problems.

\section{Main Results}
\label{sec:results}
In view of the isomorphism \eqref{eqn:iso},  Theorem \ref{thm:main} (respectively, 
Theorem \ref{thm:main-sa-closed}) gives  a 
bound (which is polynomial for fixed $d$)  on the multiplicity of the trivial representation in the $\mathfrak{S}_{\mathbf{k}}$-module $\HH^\ast(V,\F)$ (respectively, $\HH^\ast(S,\F)$).
In the current paper we generalize  both Theorems \ref{thm:main} and \ref{thm:main-sa-closed}  by proving a polynomial bound on the multiplicities of every irreducible
representation  appearing in the isotypic decomposition of $\HH^\ast(V,\F)$ and $\HH^\ast(S,\F)$ .
Note that as Example \ref{eg:basic} shows,
the dimensions of $\HH^\ast(V,\F)$, where $V$ is a 
symmetric real variety in $\R^k$ 
defined by polynomials of degree bounded by $d$ could be
exponentially large in $k$.
We also extend these basic results in several directions -- including more general actions of the symmetric group, and as a particular case symmetric varieties in $\C^k$,  as well as symmetric projective varieties.

\subsection{Affine algebraic case}
We first state our results for symmetric real algebraic subvarieties of real affine space. 
\begin{notation}
  \label{not:multiplicity-betti} 
  Let $\kk =(k_1,\ldots,k_\ell),\mm=(m_1,\ldots,m_\ell) \in \Z_{>0}^\ell$, and 
  $K = \sum_{i=1}^\ell k_i m_i$.
  For any $\mathfrak{S}_{\kk}$-symmetric
  semi-algebraic subset $S \subset \R^{K}$, any field $\F$, 
 and $\pmb{\lambda} \in \Par(\kk)$, we denote
  \begin{eqnarray*}
    m_{i,\pmb{\lambda}} (S,\F) & =  & \dim_\F \Hom_{\mathfrak{S}_\kk} (\mathbb{S}^{\pmb{\lambda}},\HH^i(S,\F))\\
    &=& \mult(\mathbb{S}^{\pmb{\lambda}}, \HH^i(S,\F)), \\
    m_{\pmb{\lambda}} (S,\F) & =  & \sum_{i}  m_{i,\pmb{\lambda}} (S,\F).
  \end{eqnarray*}
  Note that in the particular case when $\pmb{\lambda} = ((k_1),\ldots,(k_\ell))$ (i.e. when $\mathbb{S}^{\pmb{\lambda}}$ is the trivial representation of $\mathfrak{S}_\kk$), 
 \begin{eqnarray*}
    m_{i,\pmb{\lambda}} (S,\F) & =  &b_i(S/\mathfrak{S}_\kk,\F),\\
    m_{\pmb{\lambda}} (S,\F) & =  & b(S/\mathfrak{S}_\kk,\F).
  \end{eqnarray*}
  \end{notation}

\begin{notation}[Ordered tuple of degrees raised to some power]
For 
\[\dd = (d_1,\ldots,d_\ell),\m = (m_1,\ldots,m_\ell) \in \Z_{>0}^\ell,
\] 
we denote by 
\[
\dd^\m = (d_1^{m_1},\ldots,d_\ell^{m_\ell}).
\]
\end{notation}
\begin{definition}[Rank of a partition]
\label{def:rank}
For any partition $\mu \in \Par(k)$, we denote by $\rank(\mu)$ to be the length of the main diagonal in the Young diagram of 
$\mu$. Equivalently,  $\rank(\mu)$ is the side length of the largest square with a vertex at the origin  (also called the \emph{Durfee square} of $\mu$)
that fits inside the Young diagram of $\mu$ (see for example \cite[Page 65]{Stanley-vol1}).

More generally, for $\kk \in \Z_{>0}^\ell$, and $\pmb{\mu} = (\mu^{(1)},\ldots,\mu^{(\ell)}) \in \Par(\kk)$,  we define 
\[
\rank(\kk) = (\rank(\mu^{(1)}),\ldots,\rank(\mu^{(\ell)})).
\]
\end{definition}

We are now ready to state the main theorem of this section.

\begin{theorem}
\label{thm:main-product-of-symmetric-quantitative}
Let $\kk=(k_1,\ldots,k_\ell),
\m =(m_1,\ldots,m_\ell), 
\dd=(d,\ldots,d) \in \Z_{>0}^\ell$, and $K = \sum_{i=1}^{\ell} m_i k_i$.
Let 
$P \in \R[\X^{(1)},\ldots,\X^{(\ell)}]^{\mathfrak{S}_\kk}_{\leq \dd}$
be a  non-negative
polynomial
and let $V = \ZZ(P,\R^K)$.
Then, for all $\pmb{\mu} \in \Par(\kk)$,  $m_{\pmb{\mu}}(V,\F) > 0$ implies that:
\begin{enumerate}[1.]
\item
\label{itemlabel:thm:main-product-of-symmetric-quantitative1}
\[
\rank(\pmb{\mu}) \leq (2\dd)^\m, 
\]
\item
\label{itemlabel:thm:main-product-of-symmetric-quantitative2}
\begin{eqnarray*}
m_{\pmb{\mu}}(V,\F) &\leq&  \prod_{1 \leq i \leq \ell} \left(\sum_{j=0}^{(2d)^{m_i}}k_i^{O(j^2)} \binom{k_i-1}{j-1} (O(d))^{m_i j}\right)\\
&\leq& \prod_{1 \leq i \leq \ell} \left(k_i^{O((2d)^{2m_i})} (O(d))^{m_i (2d)^{m_i}}\right).
\end{eqnarray*}
In the particular case, when $\ell=1$, and $d_1=d$ and $m_1=m$  are fixed, the above bound is polynomial in $k_1=k$.
\end{enumerate}
\end{theorem}

\begin{remark}
\label{re:restriction-is-non-trivial}
Note that 
the restriction on the Specht modules that are allowed to appear in the cohomology
module $\HH^*(V,\F)$ that follows from Part  \eqref{itemlabel:thm:main-product-of-symmetric-quantitative1} 
of Theorem \ref{thm:main-product-of-symmetric-quantitative}
\emph{does not} follow only from dimension considerations,
and the Ole{\u\i}nik-Petrovski{\u\i}-Thom-Milnor bound (Theorem \ref{thm:classical}) 
on $b(V,\F)$.

For example, let $\ell=1,m_1=1,k_1=k = 2^p-1$, and let $\lambda \vdash k$ be the
partition $(2^{p-1},2^{p-2},\ldots,1)$. In this case:

\begin{eqnarray*}
\dim_\F \mathbb{S}^{\lambda} &\leq & \dim_\F \Ind_{\mathfrak{S}_\lambda}^{\mathfrak{S}_k} 
(\mbox{ since }  K(\lambda,\lambda)=1 \mbox{(Definition \ref{def:Kostka} and Proposition \ref{prop:Young}})
  \\
&=&  \binom{k}{2^{p-1},\ldots,1} \\
&\leq & O(1)^k \mbox{ using Stirling's approximation}.
\end{eqnarray*}

Thus, if $V$ is defined by a polynomial of degree bounded by $d$, and $k$ is large enough,
$\mathbb{S}^\lambda$ is not ruled out of appearing with positive multiplicity in 
$\HH^*(V,\F)$ just on the basis of the upper bound in Theorem \ref{thm:classical}.
On the other hand, it follows from
Part  \eqref{itemlabel:thm:main-product-of-symmetric-quantitative1} 
of Theorem \ref{thm:main-product-of-symmetric-quantitative} that
for all $k$ large enough,  and fixed $d$, 
\[
m_\lambda(V,\F)  = 0.
\] 
\end{remark}

We will also need the following somewhat special form of Theorem \ref{thm:main-product-of-symmetric}.

The following theorem, which yields a bound on the multiplicity of the trivial representation in $\HH^*(V,\F)$ in a special case (following the same notation as above)  follows easily from the proof of Theorem \ref{thm:main-product-of-symmetric}. It will be used later in the proof of Theorem \ref{thm:descent2-quantitative-new}.

\begin{theorem}
\label{thm:equivariant} 
Suppose that
$\mathbf{k}= (\underbrace{1, \ldots 1}_{\ell-1} ,k)$, 
$\mathbf{m}= (\underbrace{1, \ldots 1}_{\ell-1} ,m)$, 
and $(2d)^m \leq k$.
If $\pmb{\mu} = ((1),\ldots,(1),(k))$ (i.e.  $\mathbb{S}^{\pmb{\mu}}$ is the trivial representation),
then
\begin{eqnarray*}
 m_{\pmb{\mu}}(V,\F) &=&  b(V/\mathfrak{S}_\kk,\F)  \\
 &\leq & k^{(2d)^{m}} (O(d))^{m (2d)^m + \ell}.
 \end{eqnarray*}
\end{theorem}

Notice that Theorem  \ref{thm:equivariant} generalizes to the case $m > 1$, Corollary 3 in \cite{BC-advances}.

We have the following theorem for symmetric complex affine varieties.

\begin{theorem}
[Symmetric complex affine varieties]
\label{thm:main-product-of-symmetric-quantitative-complex}
Let $\kk=(k_1,\ldots,k_\ell),
\m =(m_1,\ldots,m_\ell), 
\dd=(d,\ldots,d) \in \Z_{>0}^\ell$,
and 
$K = \sum_{i=1}^{\ell} k_i m_i$.
Let 
$\mathcal{P} \subset \C[\X^{(1)},\ldots,\X^{(\ell)}]^{\mathfrak{S}_\kk}_{\leq \dd}$
be a finite set of polynomials.
Let  $V = \ZZ(\mathcal{P},\C^K)$.
Then, for all $\pmb{\mu} \in \Par(\kk)$,  $m_{\pmb{\mu}}(V,\F) > 0$ implies that:
\begin{enumerate}[1.]
\item
\[
\rank(\pmb{\mu}) \leq (4\dd)^{2\m},
\]
\item
\begin{eqnarray*}
m_{\pmb{\mu}}(V,\F) &\leq& k_i^{O(d^2)} \prod_{1 \leq i \leq \ell} \left(\sum_{j=0}^{(4d)^{2m_i}} \binom{k_i}{j} (O(d))^{2 m_i j}\right)\\
&\leq& \prod_{1 \leq i \leq \ell} \left(k_i^{O((4d)^{4m_i})} (O(d))^{2m_i (4d)^{2m_i}}\right).
\end{eqnarray*}
\end{enumerate}
\end{theorem}

\subsection{Affine semi-algebraic case}
We now state our results in the semi-algebraic case.

\begin{theorem}[Symmetric affine semi-algebraic sets]
\label{thm:main-product-of-symmetric-sa-quantitative}
Let $\kk=(k_1,\ldots,k_\ell),
\m =(m_1,\ldots,m_\ell), 
\dd=(d,\ldots,d) \in \Z_{>0}^\ell$,  $K = \sum_{i=1}^{\ell} k_i m_i$.
Let 
$\mathcal{P} \subset \R[\X^{(1)},\ldots,\X^{(\ell)}]^{\mathfrak{S}_\kk}_{\leq \dd}$
be a finite set of polynomials,
and let $\card(\mathcal{P}) = s$. Let  
$S \subset \R^K$  be a $\mathcal{P}$-closed semi-algebraic set.
Then, for all $\pmb{\mu} \in \Par(\kk)$,  $m_{\pmb{\mu}}(S,\F) > 0$ implies that:
\begin{enumerate}[1.]
\item
\[
\rank(\pmb{\mu}) \leq (4\dd)^{\m},
\]
\item
\begin{eqnarray*}
m_{\pmb{\mu}}(S,\F) &\leq& O(s)^D \prod_{1 \leq i \leq \ell} \left(k_i^{O((4d)^{2m_i})} (O(d))^{2m_i (4d)^{m_i}}\right),
\end{eqnarray*}
where 
\[
D=D (\kk,\m,d) = \sum_{i=1}^{\ell} \min (m_i k_i , d^{m_i}).
\]
\end{enumerate}

In the particular case, when $\ell=1$, and $d_1=d$ and $m_1=m$  are fixed, both bounds are polynomial in $s$ and
$k_1=k$.
\end{theorem}

\subsection{Projective case}
We can apply our results obtained in the previous section to study the topology of symmetric projective
varieties as well. 
We state one such result below.

\begin{theorem}[Symmetric complex projective varieties]
\label{thm:symmetric-complex-projective}
Let $V \subset \PP_\C^k$ be defined by a finite set of homogeneous polynomials in $\C[X_0,\ldots,X_k]^{\mathfrak{S}_{k+1}}_{\leq d}$.
Then, for all $\mu \in \Par(k+1)$,  $m_{\pmb{\mu}}(V,\F) > 0$ implies that:
\begin{enumerate}[1.]
\item
\[
\rank(\mu) \leq (4d),
\]
\item
\begin{eqnarray*}
m_{\pmb{\mu}}(S,\F) &\leq& 
k^{O(d^{4})} d^{O(d)}.
\end{eqnarray*}
\end{enumerate}
\end{theorem}

\begin{remark}
\label{rem:symmetric-complex-projective}
Suppose $V \subset \PP_\C^k$ be defined by symmetric homogeneous polynomials 
in $\C[X_0,\ldots,X_k]$ of degrees bounded by $d$.
Unlike in the affine case, it is not true that dimensions of equivariant cohomology,
$\dim_\F \HH_{\mathfrak{S}_{k+1}}^*(V,\F)$,
are bounded by a function of $d$ independent of $k$.
For example,
\[
\HH^*(\PP_\C^k,\F) \cong \HH^*(\PP_\C^k/\mathfrak{S}_{k+1},\F),
\]
 and thus
 \[\dim_\F  \HH^*(\PP_\C^k/\mathfrak{S}_{k+1},\F) = k+1,
 \]
 which clearly grows linearly with $k$.
\end{remark}

\subsection{Application to bounding topological complexity of images of polynomial maps}
\label{subsec:non-equivariant} 
In this section we discuss an application of 
Theorem  \ref{thm:equivariant}  
to bounding the Betti
numbers of images of real algebraic varieties under linear projections.
In \cite{BC-advances}, similar results were proved in the very special case of projections of the form 
$\pi:\R^{k+1} \rightarrow \R^k$. In this paper, since we consider more general actions of the symmetric group, we are able to handle projections along more than one variables, and so are
able to strengthen as well as generalize the results in \cite{BC-advances}.

In order to state our results more precisely we first introduce some notation.
Let 
$P \in \R [ Y_{1} , \ldots ,Y_{k} ,X_{1} , \ldots ,X_{m}]$ be a non-negative polynomial with 
$\deg ( P ) \leq d$.
Let $\pi : \R^{m+k}
\longrightarrow \R^{k}$ be the projection map to the first $k$ co-ordinates,
and let 
$V = \ZZ(P,\R^{m+k})$.
We
consider the problem of bounding the Betti numbers of the image 
$\pi(V)$.
Bounding the complexity of the image under projection of semi-algebraic sets is a very
important and well-studied problem related to quantifier elimination in the the first order
theory of the reals, and has many ramifications -- including in computational
complexity theory (see for example \cite{Basu-sheaf}).

There are two different approaches. One can first obtain a semi-algebraic
description of the image 
$\pi(V)$
with bounds on the degrees and the number
of polynomials appearing in this description (via results in effective quantifier elimination in the
the first order theory of the reals), and then apply known bounds on
the Betti numbers of semi-algebraic sets in terms of these parameters. Another
approach (due to Gabrielov, Vorobjov and Zell \cite{GVZ04}) 
is to use the ``descent spectral sequence'' of the map 
$\pi |_{V}$
which abuts to the cohomology of 
$\pi ( V )$
and bound the Betti numbers of
$\pi ( V )$
by bounding the dimensions of the $E_{1}$-terms of this spectral
sequence. 
For this approach it is essential that 
the map $\pi$ is proper (which is ensured by requiring that 
$V$
is 
bounded)
since in the general case the spectral sequence might not converge to $\HH^*(S,\F)$. The second approach produces a slightly better bound. The following
theorem (in the special case of algebraic sets) whose proof uses the second approach appears in {\cite{GVZ04}}.

\begin{theorem} \cite{GVZ04}
  \label{thm:descent-quantitative} 
With the same notation as above,
  \begin{eqnarray}
  \label{eqn:descent-quantitative}
    b ( \pi ( V ) ,\F ) & = & ( O ( d ) )^{( m+1 ) k} .
  \end{eqnarray}
\end{theorem}

Notice that in the exponent of the bound in  \eqref{eqn:descent-quantitative}, there is a 
factor of $(m+1)$ which is linear in the dimension of the fibers of the projection $\pi$.
This factor is also present if one uses effective quantifier elimination method to bound the Betti numbers of 
$\pi(V)$.
  Using Theorem \ref{thm:main-product-of-symmetric-sa-quantitative} we are able to remove this \emph{multiplicative} factor of $(m+1)$ in the exponent of the bound in 
\eqref{eqn:descent-quantitative}
at the expense of an extra additive term that depends just  on $d$ and $m$.

We now state the result more precisely.
In \cite{BC-advances},
the following
bound on the Betti numbers of the image under projection to a subspace of 
dimension one less than that of the ambient space of real algebraic varieties (i.e. with $m=1$), as well as of semi-algebraic sets (not
necessarily symmetric).

\begin{theorem} \cite[Theorem 10]{BC-advances}
  \label{thm:descent2-quantitative}Let $P \in \R [ Y_{1} , \ldots ,Y_{k} ,X
  ]$ be a non-negative polynomial and with $\deg ( P ) \leq d$. Let $V= \ZZ
  \left( P, \R^{k+1} \right)$ be bounded, and $\pi : \R^{k} \times \R
  \longrightarrow\R^{k}$ be the projection map to the first $k$ coordinates.
  Then,
  \begin{eqnarray*}
    b ( \pi ( V ) ,\F ) & \leq &  \left( \frac{k}{d} \right)^{2d} ( O ( d ) )^{k+2d+1}.
  \end{eqnarray*}
\end{theorem}

In this paper we generalize the above results to the case $m>1$. We prove the following theorem.
\begin{theorem}
  \label{thm:descent2-quantitative-new}
  Let $P \in \R [ Y_{1} , \ldots ,Y_{k} ,X_1,\ldots,X_m]$ be a non-negative polynomial and with $\deg ( P ) \leq d$. Let $V= \ZZ(P, \R^{k+m})$ be bounded, and $\pi : \R^{k} \times \R^{m}
  \longrightarrow\R^{k}$ be the projection map to the first $k$ coordinates.
  Then,
  \begin{eqnarray}
  \label{eqn:descent2-quantitative-new}
    b ( \pi ( V ) ,\F ) & \leq & k^{(2d)^m} (O(d))^{k+ m (2d)^m +1} .
  \end{eqnarray}
\end{theorem}

\begin{remark}
\label{rem:descent2-quantitative-new}
For every fixed $d$ and $m$ and $d,m \geq 1$,
the bound in inequality \eqref{eqn:descent2-quantitative-new} in Theorem \ref{thm:descent2-quantitative-new}  is better than the one in
\eqref{eqn:descent-quantitative} in Theorem \ref{thm:descent-quantitative},  for all large enough $k$, since
in this case 
\[
k+ m (2d)^m +1 \ll (m+1)k.
\]
\end{remark}

\subsection{Application to proving lower bounds on degrees}
The upper bounds in the theorems stated above can be potentially applied to prove lower bounds on the degrees of polynomials needed to define symmetric varieties having certain prescribed geometry.
We describe one such example.

 \begin{example}
 \label{eg:lower-bound}
 Let $k = 2^p-1$, and let 
 $\tilde{V}_k$ be any non-empty closed and bounded  semi-algebraic set contained in the subset of $\R^k$ defined by
 \begin{eqnarray*}
  X_1=&\cdots &=X_{2^{p-1}}  \\
  &\neq& \\
   X_{2^{p-1}+1} = &\cdots& = X_{2^{p-1}+2^{p-2}} \\
    &\neq&\\
  X_{2^{p-1}+2^{p-2}+1} = &\cdots& =  X_{2^{p-1}+\cdots+2^2+1} X_{2^{p-1}+\cdots+2^1} \\
   &\neq& \\
   &X_{2^{p-1}+\cdots+2^1+1}&.
  \end{eqnarray*}
  
Then, the stabilizer of $\tilde{V}_k$ under the action of $\mathfrak{S}_k$ on $\R^k$, is the
Young subgroup $\mathfrak{S}_{\lambda^{(k)}}$, where $\lambda^{(k)} = (2^{p-1},2^{p-2},\ldots,1)$.
Let $V_k$ be the orbit of $\tilde{V}_k$ under the action of $\mathfrak{S}_k$. In other words,
\[
V_k  = \mathfrak{S}_k \cdot \tilde{V}_k.
\]

Then, 
\begin{eqnarray*}
b_0(V_k,\F) & = &  b_0(\tilde{V}_k,\F) \cdot \binom{k}{2^{p-1}, 2^{p-2},\ldots,2^0}  \\
               &=&  b_0(\tilde{V}_k,\F) \cdot (\Theta(1))^k \mbox{ using Stirling's approximation}.
\end{eqnarray*}

We claim that that for any constant $d_0$, for  all $k$ large enough, $V_k$ cannot be described as the
set of real zeros of a polynomial $P \in \R[X_1,\ldots,X_k]$ with $\deg(P) \leq d_0$.
To see this observe that 
\[
\HH^0(V_k,\F) \cong_{\mathfrak{S}_k} b_0(\tilde{V}_k,\F) \cdot M^{\lambda^{(k)}} \mbox{ (cf. Definition \ref{def:Young})},
\]
with  $\lambda^{(k)}= (2^{p-1},2^{p-2},\ldots,1)$, and it follows that
$m_{0,\lambda^{(k)}}(V_k,\F) > 0$. However, 
clearly $\rank(\lambda^{(k)})$ is a strictly increasing function of $k$, and hence it follows from
Theorem \ref{thm:main-product-of-symmetric-quantitative} that
$V_k$ cannot be defined by polynomials with degrees bounded by $d_0$. 

Note that
in the case, when $\tilde{V}_k$ is a finite set of points, the same result can also be deduced from 
Proposition \ref{prop:half-degree}.
\end{example} 

\section{Preliminaries}
\label{sec:preliminaries}
Before proving the new theorems stated in the previous section we need some preliminary results.
These are described in the following subsections.

\subsection{Real closed extensions and Puiseux series}
\label{subsec:Puiseux}
In this section we recall some basic facts about real closed fields and real
closed extensions.

We will need some
properties of Puiseux series with coefficients in a real closed field. We
refer the reader to {\cite{BPRbook2}} for further details.

\begin{notation}[Field of Puiseux series]
\label{not:Puiseux}
  For $\R$ a real closed field we denote by $\R \left\langle \eps
  \right\rangle$ the real closed field of algebraic Puiseux series in $\eps$
  with coefficients in $\R$. We use the notation $\R \left\langle \eps_{1} ,
  \ldots , \eps_{m} \right\rangle$ to denote the real closed field $\R
  \left\langle \eps_{1} \right\rangle \left\langle \eps_{2} \right\rangle
  \cdots \left\langle \eps_{m} \right\rangle$. Note that in the unique
  ordering of the field $\R \left\langle \eps_{1} , \ldots , \eps_{m}
  \right\rangle$, $0< \eps_{m} \ll \eps_{m-1} \ll \cdots \ll \eps_{1} \ll 1$.
\end{notation}

\begin{notation}[Extensions]
\label{not:Ext}
  If $\R'$ is a real closed extension of a real closed field $\R$, and $S
  \subset \R^{k}$ is a semi-algebraic set defined by a first-order formula
  with coefficients in $\R$, then we will denote by $\Ext \left(S, \R'
  \right) \subset \R'^{k}$ the semi-algebraic subset of $\R'^{k}$ defined by
  the same formula. It is well-known that $\Ext \left(S, \R' \right)$ does
  not depend on the choice of the formula defining $S$ {\cite{BPRbook2}}.
\end{notation}

\begin{notation}[Balls]
  For $x \in \R^{k}$ and $r \in \R$, $r>0$, we will denote by $B_{k} (x,r)$
  the open Euclidean ball centered at $x$ of radius $r$. If $\R'$ is a real
  closed extension of the real closed field $\R$ and when the context is
  clear, we will continue to denote by $B_{k} (x,r)$ the extension $\Ext
  \left(B_{k} (x,r) , \R' \right)$. This should not cause any confusion.
\end{notation}

\subsection{Tarski-Seidenberg transfer principle}
In some proofs involving
Morse theory (see for example the proof of Lemma \ref{lem:equivariant_morseB}), where integration of gradient flows is used in an essential way, we
first restrict to the case $\R =\mathbb{R}$. After having proved the result
over $\mathbb{R}$, we use the Tarski-Seidenberg transfer theorem to extend
the result to all real closed fields. We refer the reader to 
{\cite[Chapter 2]{BPRbook2}} for an exposition of the Tarski-Seidenberg transfer
principle.

\subsection{Equivariant Morse theory}
\label{subsec:equivariant-morse-theory}

In this section we develop some basic results in equivariant Morse theory that we will need for the proof of 
Theorem \ref{thm:main-product-of-symmetric}. Since the results of this section applies to more general (finite) groups
acting on a manifold, we state and prove our results in a more general setting than what we need in this paper.
The main results from this section that will be used later are Lemma \ref{lem:equivariant_morseA} and 
Proposition \ref{prop:equivariant-morseb} which correspond to equivariant versions of the usual Morse Lemmas A and B respectively.

Let $G$ be a finite group acting on a closed and bounded semi-algebraic $S \subset \R^k$, 
defined by $Q \leq 0$, and $W = \ZZ(Q,\R^k) = \partial S$ a
bounded non-singular real algebraic hypersurface.  
Let $e:W \rightarrow \R$ be a $G$-equivariant regular 
function with isolated  non-degenerate critical points on $W$.
For each such critical point $\x$, we will denote by $\ind^-(\x)$ the dimension of the negative eigenspace
of the Hessian of $e$ at $\x$. More precisely, the Hessian $\mathrm{Hess}(e)(\x)$ is a symmetric, non-degenerate quadratic form
on the tangent space $T_p W$, and $\ind^-(\x)$ is the number of negative eigenvalues of $\mathrm{Hess}(e)(\x)$. 

Consider the set of critical points, $\mathcal{C}$, of the function $e$ restricted to $V$. 
For any subset $I \subset \R$, we will denote by $S_I = S \cap e^{-1}(I)$. 
If $I = [\infty,c]$ we will denote $S_I = S_{\leq c}$.

In the next two lemmas we will  $\R = \mathbb{R}$ since we will use properties of gradient flows.
\begin{lemma}
  \label{lem:equivariant_morseA} 
  Let $v_1<\cdots < v_N$ be the critical values of $e$ restricted to $W$.
  Then, for $1 \leq i<N$, and for each $v \in
  [ v_{i} ,v_{i+1})$, $S_{\leq v_{i}}$ is a deformation retract of $S_{\leq v}$, and the retraction
  can be chosen to be $G$-equivariant.
\end{lemma}
\begin{proof}
See for example the proof of Theorem 7.5 (Morse Lemma A) in \cite{BPRbook2} with 
$W = \ZZ(Q,\R^k)$, $a=v_{i}$ and $b=v$, noting since $W$ is symmetric and $e$ is symmetric,  the retraction of $W_{\leq v}$ to $W_{\leq v_i}$ that is constructed in the proof of
Theorem 7.5 in \cite{BPRbook2} is symmetric as well. 
\end{proof}
We also need the following equivariant version of Morse Lemma B.

\begin{lemma}
 \label{lem:equivariant_morseB}
Let $v \in e(\mathcal{C})$ be a critical value of $e$.
Let $\mathcal{C}_v^+,\mathcal{C}_v^- , \mathcal{C}_v \subset \mathcal{C}$ be defined by
\begin{eqnarray*}
\mathcal{C}_v^+ &=& \{\x \in \mathcal{C}\;\mid \; e(\x) =v, \;\langle \grad(e), \grad(Q)\rangle(\x)  >0 \}, \\
\mathcal{C}_v^- &=& \{\x \in \mathcal{C}\;\mid \; e(\x) =v, \;\langle \grad(e), \grad(Q)\rangle(\x)  <0 \}, \\
\mathcal{C}_v &=& \mathcal{C}_v^+ \;\accentset{\circ}{\cup}\;  \mathcal{C}_v^-
\end{eqnarray*}
($\accentset{\circ}{\cup}$ denotes disjoint union).
  Then for all $0 < \eps \ll 1$, $S_{\leq v+\eps}$ retracts
  $G$-equivariantly to a space 
  \[
  S_{\leq v-\eps} \cup_{B} A,
  \]
  where  
  \[(A,B) = \coprod_{\y \in \mathcal{C}_v} ( A_{\y} ,B_{\y}),
  \]
  and for each $\y \in \mathcal{C}_v$,  $(A_{\y} ,B_{\y})$ is $G$-equivariantly homotopy equivalent to the pair 
   \[(
  \mathbf{D}^{\ind^-(\y)} \times [ 0,1 ] , \partial \mathbf{D}^{\ind^-(\y)} \times [ 0,1 ] \cup
  \mathbf{D}^{\ind^-(\y)} \times \{ 1 \})
  \]
  if $\y \in \mathcal{C}_v^+$,
  or to the pair 
  \[
  (\mathbf{D}^{\ind^{-} (\y)} , \partial \mathbf{D}^{\ind^{-} (\y)}),
  \]
  if $\y \in \mathcal{C}_v^-$.
\end{lemma}
\begin{proof}
See proof of Proposition 7.19 in \cite{BPRbook2}, noting again that the retraction constructed in that proof is symmetric in case $Q$ is a symmetric polynomial and the Morse function $e$ is symmetric as well. 
\end{proof}

Now, let $\overline{\mathcal{C}}$ be a set containing a unique representative from 
each $G$-orbit of $\mathcal{C}$.
 
Let $\overline{\mathcal{C}}_i\subset \overline{\mathcal{C}}$ be the set of representatives of the different orbits corresponding to the critical value $v_i$ -- in other words, 
$e(\overline{\mathcal{C}}_i) =\{v_i\}$. 
Note that the cardinality of $\overline{\mathcal{C}}_i$ can be greater than one.
For each $\x \in \overline{\mathcal{C}}_i$, let $G_\x \subset G$ 
denote the stabilizer subgroup of $\x$.

Let also for each $i,0 \leq i \leq N$, and $0< \eps \ll 1$, $j_{i,\eps}$ denote inclusion $S_{\leq v_i -\eps} \hookrightarrow S_{\leq v_i +\eps}$, and $j_{i,\eps}^\ast: \HH^*(S_{\leq v_i +\eps},\F)  \rightarrow
\HH^*(S_{\leq v_i -\eps},\F)$ the induced homomorphism (which is in fact a homomorphism
of the corresponding $G$-modules).

Let $\overline{\mathcal{C}}_i = \{\x^{i,1},\ldots,\x^{i,N_i} \}$ (choosing an arbitrary order).

\begin{proposition}
\label{prop:equivariant-morseb}
The homomorphism $j_{i,\eps}^\ast$ factors through $N_i$ homomorphisms
as follows:
\[
\HH^*(S_{\leq v_i +\eps},\F) =M_0 \xrightarrow{j_{i,\eps,1}^\ast} M_1 \xrightarrow{j_{i,\eps,2,\ast}} \cdots
\xrightarrow{j_{i,\eps,N_i,\ast}} M_{N_i+1} =\HH^*(S_{\leq v_i -\eps},\F) ,
\]
where each $M_h$ is a finite dimensional $G$-module, and
for each $h, 1\leq h \leq N_i$,
 either
\begin{enumerate}[(a)]
\item
\label{item:theorem-main:a}
$j_{i,\eps,h}^\ast$ is injective, and 
\[
M_{h+1} \cong M_h\oplus \Ind_{G_{{\x^{i,h}}}}^{G}(W_{\x^{i,h}}),
\]
for some one-dimensional representation $W_{\x^{i,h}}$ of $G_{\x^{i,h}}$, or
\item
\label{item:theorem-main:b}
the homomorphism
$j_{i,\eps,h}^\ast$ is surjective, and 
\[
M_{h} \cong M_{h+1} \oplus \Ind_{G_{\x^{i,h}}}^{G}(W_{\x^{i,h}}),
\]
for some one-dimensional representation $W_{\x^{i,h}}$ of $G_{\x^{i,h}}$.
\end{enumerate}
\end{proposition}
 
\begin{proof}
We first assume that $\R = \mathbb{R}$.
Using Lemma \ref{lem:equivariant_morseB} (equivariant Morse Lemma B), we have that 
$S_{\leq v_i +\eps}$ can be retracted $G$-equivariantly to a 
semi-algebraic set
\[
\tilde{S}_i = S_{\leq v_i -\eps} \coprod_{\substack{1 \leq j \leq N_i,\\
 \langle\grad (e), \grad(Q)\rangle(\x^{i,j}) < 0}}
\coprod_{\y \in G \cdot \x^{i,j}} \mathbf{D}^{\ind^-(\x^{i,j})}/\sim,
\]
where the identification $\sim$ identifies the boundaries of the disks
$\mathbf{D}^{\ind^-(\x^{i,j})}$ with spheres of the same dimension in  
$S_{\leq v_i -\eps}$.

Since the different balls $\mathbf{D}^{\ind^-(\x^{i,j})}$ are disjoint, we can decompose the 
gluing process by gluing disks belonging to each orbit successively,  
choosing the order arbitrarily, and thus obtain a filtration,
\[
 S|_{e \leq v_i -\eps}=S_{i,0} \subset S_{i,1} \subset \cdots  \subset S_{i,N'_i} = \tilde{S}_i,
\]
where $S_{i,j} = S_{i,j-1} \coprod \left(\coprod_{\y \in \orbit( \x^{i,j'})}
 \mathbf{D}^{\ind^-(\y)}/\sim\right)$ for some $j', 1 \leq j' \leq N_i$.
 
 Let $D_{i,j'}$ denote the disjoint union of the balls $\mathbf{D}^{\ind^-(\x^{i,j'})}$, and 
 $C_{i,j'} \subset D_{i,j'}$ the disjoint union of their boundaries.
 
 We have the Mayer-Vietoris exact sequence
 {\small
 \[
\cdots \rightarrow \HH^{p-1}(C_{i,j'},\F) \rightarrow   \HH^{p}(S_{i,j},\F)  \rightarrow \HH^p(D_{i,j'},\F) \oplus \HH^p(S_{i,j-1},\F) \rightarrow  \HH^p(C_{i,j'},\F) \rightarrow \cdots
 \]
 }
 which is also equivariant.
 Let $n = \ind^-(\x^{i,j'})$, and assume $n \neq 0$.
 In this case, $\HH^p(C_{i,j'},\F)=0$ unless $p=0,n-1$. 
 Now, $\HH^{n-1}(C_{i,j'},\F)$ is a direct sum of $\card(\orbit(\x^{i,j'}))$, each summand
 is stable under the action of a 
 subgroup of $G$ each isomorphic to $G_{\x^{i,j'}}$, and is thus
 a one-dimensional representation  of $G_{\x^{(i,j')}}$ which we denote by
 by $W_{\x^{i,j'}}$. It follows that the representation $\HH^{n-1}(C_{i,j'},\F)$ is the induced representation $\Ind_{G_{\x^{i,j'}}}^{G}(W_{\x^{i,j'}})$.
 From the Mayer-Vietoris sequence it is evident that either
 \begin{enumerate}[(i)]
 \item
 \[
 \HH^{n}(S_{i,j},\F) = \HH^{n}(S_{i,j-1},\F) \oplus  \HH^{n-1}(C_{i,j'},\F),
 \]
 or
 \item
 \[
 \HH^{n-1}(S_{i,j},\F) \oplus \HH^{n-1}(C_{i,j'},\F) \cong \HH^{n}(S_{i,j-1},\F).
 \]
\end{enumerate}
 These two cases correspond to \eqref{item:theorem-main:a} and \eqref{item:theorem-main:b}
 respectively.
 Finally, we extend the proof to general $\R$ using the Tarski-Seidenberg transfer principle in 
 the usual way (see \cite[Chapter 7]{BPRbook2} for example).
\end{proof}

\subsection{
Structure of critical points of a symmetric Morse function on a symmetric hypersurface of small degree in $\R^k$
}
\label{subsec:degree-principle}
In this section we prove an important proposition (Proposition \ref{prop:half-degree}) that forms the basis of all our quantitative results.
It generalizes to the multi-symmetric case (i.e. for $\m$ not necessarily equal to $(1,\ldots,1)$) 
a similar result proved earlier (see \cite{Riener,Timofte03,BC-advances}).

\begin{notation}
 Let $\kk = (k_1,\ldots,k_\ell),\m = (m_1,\ldots,m_\ell),
 \mathbf{p} =(p_1,\ldots,p_\ell)  \in \Z_{>0}^\ell$, and let
 \[
 K
 = \sum_{1 \leq i \leq \ell} k_i m_i.
 \]
 
 We denote by $A_{\kk,\m}^{\mathbf{p}}$ the
  subset of $\R^{k}$ defined by
  \begin{eqnarray*}
    A^{\mathbf{p}}_{\kk,\m} & = & \left\{ \x= (\x^{(1)} , \ldots, \x^{(
    \ell)}) \mid  \ensuremath{\operatorname{card}} 
    \left(\bigcup_{j=1}^{k_i} \{ \x_{j}^{(i)} \} \right) = p_{i} \right\} .
  \end{eqnarray*}
  
  In the special case when $\ell=1,m=1$, we define for $p \leq k = k_1$,
  \begin{eqnarray*}
    A^{p}_{k} & = & \left\{ \x=(x_1,\ldots,x_k)  \mid  \ensuremath{\operatorname{card}} \left(
    \bigcup_{j=1}^{k} \{ x_j \} \right) = p \right\} .
  \end{eqnarray*}
\end{notation}

Let $\kk, 
\m, 
 \dd, 
 \in \Z_{>0}^\ell$,
  $K= \sum_{i=1}^{\ell} k_{i} m_{i}$, 
  and $P \in 
  \R [ \X^{(1)} , \ldots, \X^{(\ell)} ]^{\mathfrak{S}_\kk}_{\leq \dd}$.
 
The following proposition generalizes (to the case $\m \neq (1,\ldots,1)$) 
Proposition 5
in \cite{BC-advances}. 
\begin{proposition}
  \label{prop:half-degree}
  Let 
  \[
  e =\sum_{1 \leq h \leq \ell} \sum_{1 \leq i \leq m_h}\sum_{1 \leq j \leq k_h}
X^{(h)}_{i,j}. 
\]

Let $\mathcal{C}$ denote the set of critical points of $e$ restricted to 
$W= \ZZ(P, \R^{K})$, and suppose that $\mathcal{C}$ is a finite set. Then, 
\[
\mathcal{C} \subset \bigcup_{\mathbf{p} \leq \dd^\m} A^{\mathbf{p}}_{\kk,\m}.
\] 
\end{proposition}
\begin{proof}
For $\m= \mathbf{1} := (1,\ldots,1) $, the proposition follows immediately from 
\cite[Proposition 5]{BC-advances}.
Suppose that 
$\x = (\x^{(1)},\ldots,\x^{(\ell)}) \in \mathcal{C}$. 
For $\x = (\x^{(1)},\ldots,\x^{(\ell)}) \in \R^K$,
and $\mathbf{i} = (i_1,\ldots,i_\ell) \in [1,m_1] \times \cdots \times [1,m_\ell] $ 
denote by $\bar{\x}_{\mathbf{i}} =(\x^{(1)}_{i_1},\ldots,\x^{(\ell)}_{i_\ell}) \in \R^{K'}$,
where $K' = \sum_{i=1}^{\ell} k_i$.
It follows from the case $\mathbf{1}= (1,\ldots,1) $ 
(\cite[Proposition 5]{BC-advances}),
that 
for each $\mathbf{i} \in [1,m_1] \times \cdots \times [1,m_\ell]$,
\[
\x_{\mathbf{i}} \in \bigcup_{\mathbf{p} \leq \dd} A^{\mathbf{p}}_{\kk,\mathbf{1}}.
\]

This proves that for each $\x  = (\x^{(1)},\ldots,\x^{(\ell)}) \in \mathcal{C}$,
each row of each $(m_h \times k_h)$-matrix $\x^{(h)}$ has at most $d$ distinct entries,
and this implies that the matrix $\x^{(h)}$ has at most $d^{m_h}$ distinct columns.
This implies that  
\[
\x \in \bigcup_{\mathbf{p} \leq \dd^\m} A^{\mathbf{p}}_{\kk,\m},
\] 
which proves the proposition.
\end{proof}

\subsection{Deformation}
\label{subsec:deformation}
In this section we recall from \cite{BC-advances} an important technique for equivariantly deforming a given real variety, such that the deformed variety
has good algebraic and topological properties. The results that we are going to use later are Propositions \ref{prop:alg-to-semialg} and \ref{prop:non-degenerate} 
(both of which are reproduced here from \cite{BC-advances} for the reader's convenience)

Let $\kk =(k_1,\ldots,k_\ell),\mm=(m_1,\ldots,m_\ell) \in \Z_{>0}^\ell, K = \sum_{i=1}^\ell k_i m_i$, and $d \geq 0$. Following the notation introduced previously,
\begin{notation}[Deformation]
  \label{not:def}For any $P \in \R [ \X^{(1)} , \ldots, \X^{(\ell)} ]$ we denote
  \[ \Def (P, \zeta ,d) = P -  \zeta   \left(1+ \sum_{1\leq h\leq \ell}\sum_{1\leq i \leq m_h} \sum_{1\leq j\leq k_h}  (X^{(h)}_{i,j})^{d} \right), 
  \]
  where $\zeta$ is a new variable.
\end{notation}

Notice that if $P$ is $\mathfrak{S}_\kk$-symmetric,
then so is $\Def (P,\zeta ,d)$.

The following two propositions appear in \cite{BC-advances}. 
We restate them here for the ease of the reader.

\begin{proposition}
\cite[Proposition 3]{BC-advances}
\label{prop:alg-to-semialg}
 Let $d$ be even, 
 $P$ be $mathfrak{S}_\kk$- symmetric and non-negative polynomial,  and suppose that $V = \ZZ \left(P, \R^{K} \right)$ is bounded. 
 The variety
  $\Ext \left(V, \R \langle \zeta \rangle^{K} \right)$ is a semi-algebraic
  deformation retract of the (symmetric) semi-algebraic subset $S$ of $\R
  \langle \zeta \rangle^{K}$
  consisting
  of the union of the semi-algebraically connected components of the semi-algebraic set defined by the inequality
  \[
  \Def (P, \zeta ,d) \leq 0,
  \] 
  which are bounded over $\R$,
  and hence is semi-algebraically homotopy equivalent to $S$.
 \end{proposition}

\begin{proposition}
\cite[Proposition 4]{BC-advances}
 \label{prop:non-degenerate}
  Let 
  $P \in \R [ \X^{(1)} , \ldots, \X^{(\ell)} ]$,
  and
  $d$ be an even number with $\deg (P) 
  < 
  d=p+1$, with $p$ a prime. Let
  \[
  e = \sum_{1\leq h\leq \ell}\sum_{1 \leq i\leq m_h} \sum_{1\leq j\leq k_h}  
  X^{(h)}_{i,j},
  \] and
  \[
  V_{\zeta} = \ZZ \left(\Def (P,\zeta ,d) , \R \langle \zeta \rangle^{K} \right).
  \] 
  Suppose also that $\gcd(p,K) =1$. Then, the critical points of $e$ restricted to $V_{\zeta}$ are
  finite in number, and each critical point is non-degenerate.
\end{proposition}

\subsection{Representation theory of products of symmetric groups}
\label{subsec:representation-of-Sn}
In this section we recall some well known facts from the representation theory of symmetric groups,  and prove one new result 
(Proposition \ref{prop:multiplicity})
that will be used later.
The following classical formula (due to Frobenius) gives the dimensions of the representations $\mathbb{S}^{\lambda}$ in terms of the \emph{hook lengths} of the partition $\lambda$ defined below.

\begin{definition}[Hook lengths]
Let $B(\lambda)$ denote the set of boxes in the  Young diagram corresponding to a partition 
$\lambda \vdash k$. For a box $b \in B(\lambda)$, the length of the hook of $b$, denoted 
$h_b$ is the number of boxes strictly to the right and below $b$ plus $1$.
\end{definition}

\begin{theorem}[Hook length formula]
\label{thm:hook}
Let $\lambda \vdash k$. Then,
\begin{eqnarray}
\label{eqn:hook}
\dim_\F \mathbb{S}^{\lambda} &=&  \frac{k!}{\prod_{b \in B(\lambda)} h_b}.
\end{eqnarray}
\end{theorem}

\begin{definition}[Young module]
\label{def:Young}
For $\lambda \vdash k$, we will denote 
\[
M^\lambda = \Ind_{\mathfrak{S}_\lambda}^{\mathfrak{S}_k} (\mathbf{1}_{\mathfrak{S}_\lambda})
\] 
(where $\mathbf{1}_{\mathfrak{S}_\lambda}$ denotes the trivial one-dimensional representation
of $\mathfrak{S}_\lambda$).
\end{definition}

\begin{definition}[Dominance order]
\label{def:dominance}
For any two partitions $\mu=(\mu_1,\mu_2,\ldots),\lambda=(\lambda_1,\lambda_2,\ldots) \in \Par(k)$,
we say that $\mu \gdom \lambda$, if for each $i\geq 0$, $\mu_1+\cdots+\mu_i \geq \lambda_1+\cdots+\lambda_i$.
This is a partial order on $\Par(k)$.
More generally, for $\kk = (k_1,\ldots,k_\ell) \in \Z_{>0}^\ell$,
and $\pmb{\mu} = (\mu^{(1)},\ldots,\mu^{(\ell)}), \pmb{\lambda} = (\lambda^{(1)},\ldots,\lambda^{(\ell)}) \in \Par(\kk)$, we denote
$\pmb{\mu} \gdom \pmb{\lambda}  $ if and only if $ \mu^{(i)} \gdom \lambda^{(i)} $ for each $i, 1\leq i \leq \ell$.
\end{definition}

We also need the definitions of Kostka numbers and the Littlewood-Richardson coefficients.

\begin{definition}[Kostka numbers]
\label{def:Kostka}
For $\lambda,\mu \vdash k$, $K(\mu,\lambda)$ denotes the number of semi-standard Young tableaux of
shape $\mu$ and weight $\lambda$ (see \cite{Procesi-book} for definitions of semi-standard Young tableaux, and 
also their shape and weight). 
\end{definition}

The following fact is very basic  (see for example \cite[Theorem 3.6.11]{Ceccherini-book} or \cite[page 541, \S 7.3]{Procesi-book}).

\begin{proposition}[Young's rule]
\label{prop:Young}
Let $k \in \N$, and $\lambda\in \Par(k)$. Then,
\[
\Ind_{\mathfrak{S}_\lambda}^{\mathfrak{S}_k} \big(\mathbb{S}^{(\lambda_1)}\boxtimes \cdots \boxtimes \mathbb{S}^{(\lambda_{\length(\lambda)})} \big) \cong \bigoplus_{\mu \gdom  \lambda} K(\mu,\lambda) \mathbb{S}^{\mu}.
\]
\end{proposition}

\begin{definition}[Littlewood-Richardson coefficients]
\label{def:LR}
For $\lambda \vdash m, \mu \vdash n, \nu \vdash m+n$, $c^\nu_{\lambda,\mu}$ is the multiplicity
of the irreducible representation $\mathbb{S}^\nu$ in 
$\Ind_{\mathfrak{S}_m \times \mathfrak{S}_n}^{\mathfrak{S}_{m+n}}(\mathbb{S}^\lambda \boxtimes \mathbb{S}^\mu)$.
\end{definition}

In order to state the main new result in this section  (Proposition \ref{prop:multiplicity} below) we need one more notation.

\begin{notation}[Induced representations and multiplicities]
\label{not:induced-rep-and-mult}
For each $\lambda \vdash k$,
we denote by $\overline{\Par}(\lambda)$ the set of partitions $\mu \vdash k$ such that,
there exists a decomposition $\lambda = \lambda' \coprod \lambda''$, 
$\lambda' = (\lambda'_1,\ldots,\lambda'_{\ell'}), \lambda''=(\lambda''_1,\ldots,\lambda''_{\ell''}), \ell'+\ell''=\ell = \length(\lambda)$,such that $\mathbb{S}^\mu$ occurs with positive multiplicity
in the representation 
\[
\mathbb{S}_{\lambda',\lambda''} := 
\Ind_{\mathfrak{S}_{\lambda'} \times \mathfrak{S}_{\lambda''}}^{\mathfrak{S}_k}\left(\left(\boxtimes_{i=1}^{\ell'} \mathbb{S}^{(\lambda'_i)}\right) \boxtimes \left(\boxtimes_{j=1}^{\ell''} \mathbb{S}^{(1^{\lambda''_j})}\right) \right),
\] 
and we denote the multiplicity of $\mathbb{S}^\mu$ in $\mathbb{S}_{\lambda',\lambda''}$  by $m^\mu_{\lambda',\lambda''}$.

More generally, for $\kk \in \Z_{>0}^\ell$, and $\pmb{\lambda}= (\lambda^{(1)},\ldots,\lambda^{(\ell)}) \in \Par(\kk)$, we denote by 
$\overline{\Par}(\pmb{\lambda}) = \overline{\Par}(\lambda^{(1)}) \times \cdots \times \overline{\Par}(\lambda^{(\ell)})$.
\end{notation}

\begin{proposition}
\label{prop:multiplicity}
Let $k,d > 0$, $\lambda \in \Par(k,d)$
such that $\lambda = \lambda' \coprod \lambda''$,
and $\mu  \in \overline{\Par}(\lambda)$.
Then,
\begin{enumerate}[1.]
\item
\label{item:prop:multiplicity1}
\[
\rank(\mu) \leq d,
\]
\item
\label{item:prop:multiplicity2}
\begin{eqnarray}
\label{eqn:prop:multiplicity2}
m^\mu_{\lambda',\lambda''} &=& 
\sum_{ \substack{\nu' \vdash |\lambda'|,\nu' \gdom \lambda' \\ \nu'' \vdash |\lambda''|,\nu'' \gdom \widetilde{\lambda''}}} K(\nu',\lambda')\cdot K(\nu'',\widetilde{\lambda''})\cdot c^{\mu}_{\nu',\nu''},
\end{eqnarray}
\item
\label{item:prop:multiplicity3}
\begin{eqnarray*}
\label{eqn:prop:multiplicity3}
\sum_{\mu \vdash k} m^\mu_{\lambda',\lambda''}&\leq& 
k^{O(d^2)}.
\end{eqnarray*}
\end{enumerate}
\end{proposition}

\begin{remark}
It is well known that
$K(\mu,\mu) = 1$ for all $\mu \in \Par(k)$, $K(\mu,\lambda)= 0$ unless $\mu \gdom \lambda$. Finally, if $\mu$ is the maximal element in the dominance ordering $\gdom$ on $\Par(k)$, that is $\mu = (k)$, then $K(\mu,\lambda) = 1$ for all $\lambda \in \Par(k)$.
In particular, in conjunction with Schur's lemma the above fact implies, that the trivial representation, $\mathbb{S}^{(k)}$ occurs with multiplicity equal to $1 (= K((k),\lambda))$ in $\Ind_{\mathfrak{S}_\lambda}^{\mathfrak{S}_k}\big(\boxtimes_{j=1}^{\length(\lambda)}\mathbb{S}^{(\lambda_j)}\big)$.
\end{remark}

\begin{remark}
Note also that the representation
$\Ind_{\mathfrak{S}_\lambda}^{\mathfrak{S}_k}(\boxtimes_{j=1}^{\length(\lambda)} \mathbb{S}^{(\lambda_j)})$
is isomorphic to the permutation representation
of $\mathfrak{S}_k$ on the set of cosets $\mathfrak{S}_k/\mathfrak{S}_\lambda$, and in particular
\[
\dim_\F \Ind_{\mathfrak{S}_\lambda}^{\mathfrak{S}_k}(\boxtimes_{j=1}^{\length(\lambda)} \mathbb{S}^{(\lambda_j)})
=
\frac{k!}{\prod_{1\leq j \leq \length(\lambda)} \lambda_j!}.
\]
\end{remark}

In order to prove Proposition \ref{prop:multiplicity} we need  the following definition and results which are all well known.

\begin{definition}[Skew partitions, horizontal and vertical strips]
\label{def:skew-partition-strips}
For any two partitions, $\lambda=(\lambda_1,\lambda_2,\ldots) \vdash m$, $\mu=(\mu_1,\mu_2,\ldots) \vdash n$, $m \leq n$,  we say that $\lambda \subset \mu$,
if $\lambda_i \leq \mu_i$ for all $i$.

Identifying $\lambda,\mu$ with their respective Young diagrams,
we say that the \emph{skew partition} $\mu/\lambda$ is a \emph{horizontal strip} if no two cells
of $\mu/\lambda$ belong to the same \emph{column}. 
We say that $\mu/\lambda$ is a \emph{vertical strip}
if no two cells of $\mu/\lambda$ belong to the same \emph{row}.
\end{definition}

\begin{proposition}[Pieri's rule]
\label{prop:Pieri}
For $\lambda \vdash m$, and $n \geq 0$, we have the two following relations.
\begin{eqnarray*}
\Ind_{\mathfrak{S}_m \times \mathfrak{S}_n}^{\mathfrak{S}_{m+n}} (\mathbb{S}^\lambda \boxtimes \mathbb{S}^{(n)}) 
&\cong& \bigoplus_{\substack{\mu \vdash m+n \\ \mu/\lambda \mbox{ is a horizontal strip}}} \mathbb{S}^\mu, \\
\Ind_{\mathfrak{S}_m \times \mathfrak{S}_n}^{\mathfrak{S}_{m+n}} (\mathbb{S}^\lambda \boxtimes \mathbb{S}^{1^n}) 
&\cong& \bigoplus_{\substack{\mu \vdash m+n \\ \mu/\lambda \mbox{ is a vertical strip}}} \mathbb{S}^\mu.
\end{eqnarray*}
\end{proposition}

We also have the following associativity relationship that allows us to apply Pieri's rule (Proposition \ref{prop:Pieri})
iteratively.

\begin{proposition}
\label{prop:associative}
Let $n = m_1+\cdots+m_\ell$, where for each $i, 1\leq i \leq \ell$. Then,
\[
\Ind_{\mathfrak{S}_{m_1} \times\cdots \times \mathfrak{S}_{m_\ell}}^{\mathfrak{S}_{n}} 
(V_1 \boxtimes \cdots \boxtimes V_\ell)
\]
is isomorphic to 
\[
\Ind_{\mathfrak{S}_{m_1+\ldots+m_{\ell-1}} \times \mathfrak{S}_{m_\ell}}^{\mathfrak{S}_n}
(\Ind_{\mathfrak{S}_{m_1} \times \cdots \times \mathfrak{S}_{m_{\ell-1}}}^{\mathfrak{S}_{m_1+\cdots+m_{\ell-1}}}
(V_1 \boxtimes \cdots \boxtimes V_{\ell-1}) \boxtimes V_\ell),
\]
where for each $i, 1\leq i \leq \ell$, $V_i$ is an $\mathfrak{S}_{m_i}$-module.
\end{proposition}

\begin{proof}[Proof of Proposition \ref{prop:multiplicity}]
We first prove \eqref{item:prop:multiplicity2}.
Let $k'=|\lambda'|$ and $k'' = |\lambda''|$. Then, using 
Young's rule (Proposition \ref{prop:Young})
\begin{eqnarray*}
\Ind_{\mathfrak{S}_{\lambda'}}^{\mathfrak{S}_{k'}}\left(\boxtimes_{i=1}^{\ell'} \mathbb{S}^{(\lambda'_i)}\right)
&\cong&
\bigoplus_{\nu' \vdash k',\nu' \gdom \lambda'} K(\nu',\lambda')\mathbb{S}^{\nu'}, \\
\Ind_{\mathfrak{S}_{\lambda''}}^{\mathfrak{S}_{k''}}\left(\boxtimes_{i=1}^{\ell''} \mathbb{S}^{1^{\lambda''_i}}\right)
&\cong&
\bigoplus_{\nu'' \vdash k',\nu'' \gdom \widetilde{\lambda''}} K(\nu'',\widetilde{\lambda''})\mathbb{S}^{\nu''}.
\end{eqnarray*}
It follows that
\begin{eqnarray*}
 \Ind_{\mathfrak{S}_{\lambda'} \times \mathfrak{S}_{\lambda''}}^{\mathfrak{S}_{k'} \times \mathfrak{S}_{k''}}\left(\left(\boxtimes_{i=1}^{\ell'} \mathbb{S}^{(\lambda'_i)}\right) \boxtimes \left(\boxtimes_{j=1}^{\ell''} \mathbb{S}^{(1^{\lambda''_j})}\right) \right)
 \end{eqnarray*}
 is isomorphic to 
 \begin{eqnarray*}
 \bigoplus_{ \substack{\nu' \vdash k', \nu' \gdom \lambda' \\ \nu'' \vdash k'' , \nu'' \gdom \widetilde{\lambda''}}} K(\nu',\lambda')K(\nu'',\widetilde{\lambda''}) \mathbb{S}^{\nu'} \boxtimes \mathbb{S}^{\nu''}.
 \end{eqnarray*}
Eqn. \eqref{eqn:prop:multiplicity2} then follows from the isomorphism
\begin{eqnarray*}
\mathbb{S}_{\lambda',\lambda''} &\cong&
\Ind_{\mathfrak{S}_{k'} \times \mathfrak{S}_{k''}}^{\mathfrak{S}_k} \Ind_{\mathfrak{S}_{\lambda'} \times \mathfrak{S}_{\lambda''}}^{\mathfrak{S}_{k'} \times \mathfrak{S}_{k''}}\left(\left(\boxtimes_{i=1}^{\ell'} \mathbb{S}^{(\lambda'_i)}\right) \boxtimes \left(\boxtimes_{j=1}^{\ell''} \mathbb{S}^{(1^{\lambda''_j})}\right) \right)
\end{eqnarray*}
and the definition of the Littlewood-Richardson's coefficients, $c^{\mu}_{\nu',\nu''}$ (Definition \ref{def:LR}).

An alternative way of obtaining the multiplicities $m^{\mu}_{\lambda',\lambda''}$ is by applying
Pieri's rule iteratively at most  $d$ times using Propositions \ref{prop:associative} and \ref{prop:Pieri}.
Let $\length(\lambda') = \ell',\length(\lambda'') = \ell''$, so that 
$\ell' + \ell'' = \length(\lambda) \leq d$. 

Let for $1 \leq i \leq  \ell'$,
\[
M_i = \Ind_{\mathfrak{S}_{\lambda'_1 + \ldots +\lambda'_{i-1}} \times \mathfrak{S}_{\lambda'_i}}^{\mathfrak{S}_{\lambda'_1 + \ldots +\lambda'_{i}}}(M_{i-1}\boxtimes \mathbb{S}^{(\lambda'_i)}),
\]
with the convention that $M_0 = \mathbf{1}$.
For $\nu \vdash \lambda'_1 + \ldots +\lambda'_{i}$, let $m^\nu_i$ denote the multiplicity of 
$\mathbb{S}^\nu$ in $M_i$, and $m_i = \sum_{\nu \vdash \lambda'_1 + \ldots +\lambda'_{i}}  m^\nu_i$.
We prove by induction on $i$ the following two statements.
\begin{enumerate}[(a)]
\item
\label{item:multiplicity:a}
\[
m_i \leq m_{i-1} \cdot \binom{\lambda'_i+i-1}{i-1}.
\]
\item
\label{item:multiplicity:b}
For each $\nu \vdash \lambda'_1 + \ldots +\lambda'_{i}$, such that
$m^\nu_i > 0$, $\length(\nu) \leq i$.
\end{enumerate}

Assuming statements \eqref{item:multiplicity:a} and \eqref{item:multiplicity:b} hold for $i-1$ we prove them for $i$.
By induction for each $\nu' \vdash \lambda'_1 + \ldots +\lambda'_{i-1}$ with $m^{\nu'}_{i-1} > 0$,
$\length(\nu') \leq i-1$. Applying Pieri's rule (Proposition \ref{prop:Pieri}) we obtain 
\begin{eqnarray}
\label{eqn:multiplicity:1}
\Ind_{\mathfrak{S}_{\lambda'_1+\cdots+\lambda'_{i-1}}\times \mathfrak{S}_{\lambda'_{i}}}^{\mathfrak{S}_{\lambda'_1+\cdots+\lambda'_i}} (\mathbb{S}^{\nu'} \boxtimes \mathbb{S}^{(\lambda'_i)}) 
&\cong& \bigoplus_{\substack{\nu \vdash \lambda'_1+\cdots+\lambda'_i \\ \nu/\nu' \mbox{ is a horizontal strip}}} \mathbb{S}^\nu.
\end{eqnarray}
Observe that each choice of $\nu \vdash \lambda'_1+\cdots+\lambda'_i$ such that  
$\nu/\nu'$is a horizontal strip, corresponds uniquely to a composition of $\lambda'_i$ into at most
$\length(\nu')$ parts, and the number of such compositions is clearly bounded by 
\[
\binom{\lambda'_i+\length(\nu')-1}{\length(\nu')-1} \leq \binom{\lambda'_i+i-1}{i-1},
\]
since $\length(\nu') \leq i-1$ by induction hypothesis.

This proves part \eqref{item:multiplicity:a}. 
Part \eqref{item:multiplicity:b} also follows from \eqref{eqn:multiplicity:1} noting that the length of each $\mu$ that occurs on the 
right is at most $\length(\mu')+1$ which is $\leq i$ using the induction hypothesis.
This completes the proof of parts \eqref{item:multiplicity:a}  and \eqref{item:multiplicity:b}.
 
Now let for $1 \leq j \leq  \ell''$,
\[
N_j = \Ind_{\mathfrak{S}_{\ell' + \lambda''_1 + \ldots +\lambda''_{i-1}} \times \mathfrak{S}_{\lambda''_j}}^{\mathfrak{S}_{\ell' + \lambda''_1 + \cdots +\lambda''_{j}}}(N_{j-1}\boxtimes \mathbb{S}^{(\lambda''_j)}),
\]
with the convention that $N_0 = M_{\ell'}$.
For $\nu \vdash \lambda''_1 + \ldots +\lambda''_{j}$, let $n^\nu_j$ denote the multiplicity of 
$\mathbb{S}^\nu$ in $N_j$, and $n_j = \sum_{\nu \vdash \lambda''_1 + \cdots +\lambda''_{j}}  n^\nu_j$.

The following two statements are easily proved using induction on $j$. The proofs are very similar to the proofs
of \eqref{item:multiplicity:a} and \eqref{item:multiplicity:b} above and are omitted.

\begin{enumerate}[(a)]
\addtocounter{enumi}{2}
\item
\label{item:multiplicity:c}
\[
n_j \leq n_{j-1} \cdot \binom{\lambda''_j+\ell'+j-1}{\ell'+j-1}.
\]
\item
\label{item:multiplicity:d}
For each $\nu \vdash \lambda''_1 + \cdots +\lambda''_{j}$, such that
$n^\nu_j > 0$, $\length(\widetilde{\nu}) \leq \ell'+j$.
\end{enumerate}

It follows from \eqref{item:multiplicity:a}, \eqref{item:multiplicity:b}, \eqref{item:multiplicity:c}, and \eqref{item:multiplicity:d}, that
\[
\sum_{\mu\vdash k} m^\mu_{\lambda',\lambda''} \leq k^{O(d^2)},
\]
which proves \eqref{item:prop:multiplicity3}.
Finally, it is easy to check that for each $\mu$ with $m^{\mu}_{\lambda',\lambda''} > 0$
that arises in the above process satisfies
\[
\card(\{ i \mid \mu_i \geq d \}) \leq d, \\
\card(\{j \mid \tilde{\mu}_j \geq d\}) \leq d,
\]
which proves \eqref{item:prop:multiplicity1}.
\end{proof} 

\begin{remark}
\label{rem:multiplicity}
The following particular case of Proposition \ref{prop:multiplicity} will be of interest.
If $\mu=(k)$, 
\[m^\mu_{\lambda',\lambda''} = 1.
\]
\end{remark}

\subsection{Equivariant Poincar\'e duality}
\label{subsec:Poincare-duality}

In this section, we derive an equivariant version of Poincar\'e duality for oriented manifolds (Theorem \ref{thm:poincare-duality}) that was used in analyzing 
Example \ref{eg:basic}.

\begin{theorem}
\label{thm:poincare-duality}
Let $V \subset \R^k$ be a closed and bounded non-singular  semi-algebraic oriented hypersurface, which is 
stable under the standard action of $\mathfrak{S}_k$ on $\R^k$. Then, for each $p,0 \leq p \leq k$,
there is an $\mathfrak{S}_k$-module isomorphism
\[
\HH^p(V,\F) \xrightarrow{\sim} \HH^{k-p-1}(V,\F) \otimes 
\mathbf{sign}_k.
\] 
\end{theorem}

\begin{proof}
If $M$ is a $C^0$-manifold of dimension $\ell$, then the following sheaf-theoretic statement
of Poincar\'e duality is well known (see for example \cite[Corollary 5.5.6]{Schapira-notes}).

\begin{equation}
\label{eqn:iso1}
\Hom_\F(\HH^*_c(M;\F_M),\F)  \cong \HH^*(M;\mathrm{or}_M)[\ell].
\end{equation}

In our case, with $M=V$.  The $\mathfrak{S}_k$-action on the ambient space $\R^k$, induces
an $\mathfrak{S}_k$-module structure on $\HH^*(V;\F_V)$ by the induced isomorphisms
$\pi^*: \HH^*(V;\F_V) \xrightarrow{\sim} \HH^*(V;\F_V), \pi \in \mathfrak{S}_k$.

Now for $\pi \in \mathfrak{S}_k$ (and also denoting by $\pi$ the induced map $\pi:V \rightarrow V$), we have that $\pi$ induces the sign representation 
on the  
one dimensional vector space, $\Gamma(V; \mathrm{or}_V)$,  of global sections of the orientation sheaf on $V$.  
This implies the following $\mathfrak{S}_k$-isomorphism for each $p \geq 0$,
\begin{equation}
\label{eqn:iso2}
\HH^p(V;\mathrm{or}_V) \cong \HH^p(V;\F_V) \otimes 
\textbf{sign}_k.
\end{equation}
The theorem follows from \eqref{eqn:iso1} and \eqref{eqn:iso2}, after noting that since $V$ is assumed to be closed and bounded
\[
\Hom_\F(\HH^*(V,\C),\F) \cong \HH^*_c(V,\F) \cong  \HH^*(V,\F) 
\]
where all isomorphisms are $\mathfrak{S}_k$-module isomorphisms.
\end{proof}

\subsection{Equivariant Mayer-Vietoris inequalities}
In this section we derive equivariant versions of Mayer-Vietoris  inequalities that we will need to obtain bounds on the multiplicities of the various Specht-modules in the
cohomology modules of symmetric varieties and semi-algebraic sets that we consider.  We will use Propositions \ref{7:prop:prop1} and \ref{prop:MV}.

Suppose that $S_1,S_2 \subset \R^K$ are $\mathfrak{S}_{\kk}$-symmetric closed semi-algebraic sets.
Then $S_1 \cup S_2$, and $S_1 \cap S_2$ are also $\mathfrak{S}_{\kk}$-symmetric closed semi-algebraic sets, and  there is the classical Mayer-Vietoris exact sequence,
{\small
\[
\cdots \rightarrow \HH^{i}(S_1\cup S_2,\F) \rightarrow \HH^{i}(S_1,\F) \oplus \HH^{i}(S_2,\F) \rightarrow \HH^i(S_1\cap S_2,\F) \rightarrow \HH^{i+1}(S_1 \cup S_2,\F) \rightarrow \cdots 
\]  
}
where all the homomorphisms are $\mathfrak{S}_{\kk}$-equivariant. Denoting by $\HH^*(S,\F)_{\pmb{\mu}}$ the isotypic component of $\HH^*(S,\F)$ corresponding to $\pmb{\mu} \in \Par(\kk)$ for any $\mathfrak{S}_{\kk}$-symmetric closed semi-algebraic set $S \subset\R^k$, we
obtain using Schur's lemma for each $\pmb{\mu} \in \Par(\kk)$, an exact sequence,
{\small
\[
\cdots \rightarrow \HH^{i}(S_1\cup S_2,\F)_{\pmb{\mu}} \rightarrow \HH^{i}(S_1,\F)_{\pmb{\mu}} \oplus \HH^{i}(S_2,\F)_{\pmb{\mu}} \rightarrow \HH^i(S_1\cap S_2,\F)_{\pmb{\mu}} \rightarrow \HH^{i+1}(S_1 \cup S_2,\F)_{\pmb{\mu}} \rightarrow \cdots 
\]  
}

The following inequalities follow from the above exact sequence (the proofs are similar to the 
non-equivariant case and can be found in \cite{BPRbook2}).

Let $S_{1} , \ldots ,S_{s} \subset \R^{K}$, $s \ge 1$, be $\mathfrak{S}_{\kk}$-symmetric closed
semi-algebraic sets of $\R^{K}$, contained in a $\mathfrak{S}_{\kk}$-symmetric 
closed semi-algebraic set $T$.

For $1 \leq t \leq s$, let $S_{\le t} = \bigcap_{1 \leq j \leq t} S_{j}$, and
$S^{\le t} = \bigcup_{1 \leq j \leq t} S_{j}$.
Also, for $J \subset \{1, \ldots ,s\}$, $J \neq \emptyset$, let $S_{J} =
\bigcap_{j \in J} S_{j}$, and
$S^{J} = \bigcup_{j \in J} S_{j}$. Finally, let $S^{\emptyset} =T$.

\begin{proposition}
  \label{7:prop:prop1}
   \begin{enumerate}[(a)]
    \item 
    \label{7:prop:prop1:itema}
    For $\pmb{\mu} \in \Par(\kk)$ and $i \geq 0$,
    
     \begin{equation*}
      m_{i,\pmb{\mu}} (S^{\le s} ,\F) \leq \sum_{j=1}^{i+1}
      \sum_{\substack{
        J \subset \{ 1, \ldots ,s \}\\
        \card (J) =j}} 
       m_{i-j+1,\pmb{\mu}} (S_{J} ,\F) .
    \end{equation*}
    
    \item 
    \label{7:prop:prop1:itemb}
    For $\pmb{\mu} \in \Par(\kk)$ and $0 \le i \le K$,
    
     \begin{equation*}
      \label{7:eqn:prop1} m_{i,\pmb{\mu}} (S_{\le s} ,\F) \leq \sum_{j=1}^{K-i}
      \sum_{\substack{
        J \subset \{ 1, \ldots ,s \}\\
        \card (J) =j}}
      m_{i+j-1,\pmb{\mu}} (S^{J} ,\F) + \binom{s}{K-i} m_{K,\pmb{\mu}}
      (S^{\emptyset} ,\F) .
    \end{equation*}
  \end{enumerate}
\end{proposition}

\begin{proof} Follows from the proof of 
\cite[Proposition 7.33]{BPRbook2} and Schur's lemma.
\end{proof}

\begin{proposition}
  \label{prop:MV}
  If $S_{1} ,S_{2}$ are $\mathfrak{S}_{\kk}$-symmetric  closed semi-algebraic sets, then for
  $\pmb{\mu} \in \Par(\kk)$, any field $\F$ and every $i \geq 0$
  \begin{eqnarray*}
    m_{i,\pmb{\mu}} (S_{1} ,\F) +m_{i,\pmb{\mu}} (S_{2} ,\F) & \leq & m_{i,\pmb{\mu}}
    (S_{1} \cup S_{2} ,\F) + m_{i,\pmb{\mu}} (S_{1} \cap S_{2} ,\F).
    \end{eqnarray*}
\end{proposition}
\begin{proof}
It follows from the proof of 
\cite[Proposition 6.44]{BPRbook2} and Schur's lemma.
\end{proof}

\subsection{Descent spectral sequence}
\label{subsec:descent}
In this section we derive an improvement on an inequality first obtained in \cite{GVZ04} by taking advantage of the symmetry of the 
fibered products. The main result  of this section that we will use later is Theorem \ref{thm:descent2} (which will be used in the proof of 
Theorem \ref{thm:descent2-quantitative-new}).

Suppose that $V \subset \R^{k+m}$ is a closed and bounded semi-algebraic set, and $\pi: V\rightarrow Y = \pi(V)$ is the projection on the first $k$ coordinates restricted to $V$. Following \cite{GVZ04} we define for each $p \geq 0$,
\begin{equation}
\label{eqn:definition-of-fibered-product}
W_\pi^{(p)}(V) = \underbrace{X \times_\pi
    \cdots \times_\pi X}_{p+1} = \{ (\y,\x_0,\ldots,\x_p) \in \R^{k+(p+1)m}  \mid (\y,\x_i) \in V, 0 \leq i \leq p \}.
\end{equation}

Notice that $W_\pi^{(p)}(V)$ is $\mathfrak{S}_{\kk(p)}$-symmetric semi-algebraic set, where 
\[
\kk(p) = (\underbrace{1,\ldots,1}_{k},p+1),
\] and 
$\mathfrak{S}_{\kk(p)}$ acts by permuting the blocks $\x_0,\ldots,\x_p$, and by identity on the remaining coordinates.

We have the following theorem.
Following the same notation as above:
\begin{theorem}
  \label{thm:descent2} 
  \begin{eqnarray*}
    b (\pi (V) ,\F) & \leq & 
        \sum_{0 \leq p<k} b_{\mathfrak{S}_{\kk(p)}} (W_\pi^{(p)}(V) ,\F).
\end{eqnarray*}
\end{theorem}

\begin{proof}
It is proved in \cite{GVZ04}, that there exists a spectral sequence whose $E^1$ term is given by
$E^1_{p,q} = \HH_q( W_\pi^{(p)}(X),\F)$, 
and such that it converges to $\HH_{p+q}(Y,\F)$ in  a finite number of steps.
Note that  from the definition of the spectral sequence it is easy to see that $E^1_{p,q}$ has the structure of an $\mathfrak{S}_{p+1,1^k}$-module, and so does each $E^r_{p,q}$ which are all sub-quotients of $E^1_{p,q}$.
There is an isomorphism,
\[
F_n: \bigoplus_{p+q = n} E^\infty_{p,q} \rightarrow \HH_{n}(Y,\F),
\]
and $F_n$ restricted to $E^\infty_{p,q}$ is an $\mathfrak{S}_{\kk(p)}$-module isomorphism onto its image, where the $\mathfrak{S}_{\kk(p)}$-module
structure on the image is the trivial one. This implies by Schur's lemma that 
\[
E^\infty_{p,q} = (E^\infty_{p,q})^{\mathfrak{S}_{\kk(p)}}, 
\]
and also that 
\[
\dim_\F(E^\infty_{p,q}) \leq  \dim_\F(E^1_{p,q})^{\mathfrak{S}_{\kk(p)}}.
\]
Finally, observe that 
\[
(E^1_{p,q})^{\mathfrak{S}_{\kk(p)}} \cong   
\HH_q(W^{(p)}_\pi(X),\F)^{\mathfrak{S}_{\kk(p)}}.
\]
The theorem follows after observing that 
\[
\dim_\F(\HH_q(W^{(p)}_\pi(X),\F)^{\mathfrak{S}_{\kk(p)}}) = b_{\mathfrak{S}_{\kk(p)}}(W^{(p)}_\pi(X),\F).
\]
\end{proof}

\section{Proofs of the main theorems}
\label{sec:proofs-of-main}

We first prove a structural result that will be used in the proofs of the main theorem.
\subsection{Structural result}
We define for  $\kk,\dd ,\m\in \Z_{>0}^\ell$, 
a subset  $\mathcal{I}(\kk,\dd,\mm) \subset \Par(\kk)$, having the property that,
only the irreducible representations of $\mathfrak{S}_\kk$ associated to the elements from $\mathcal{I}(\kk,\dd,\mm)$
can appear in the cohomology modules of symmetric varieties in $\R^K$ defined by a non-negative 
polynomial having degree bounded by $\dd$.

\begin{definition}[Definition of $\mathcal{I}(\kk,\dd,\mm)$]
\label{def:set-of-irred}
For $\kk \in \Z_{>0}^\ell, \pmb{\lambda} \in \Par(\kk)$, we denote 
(cf. Notation \ref{not:induced-rep-and-mult})

\begin{equation}
\label{eqn:set-of-irred-lambda}
\mathcal{I}(\pmb{\lambda}) = 
\bigcup_{\substack{\lambda^{(i)} = {\lambda^{(i)}}' \coprod {\lambda^{(i)}}''\\ 1 \leq i \leq \ell}}
\{
 {\pmb{\mu}} \mid \pmb{\mu} =(\mu^{(1)},\ldots,\mu^{(\ell)})  \in \Par(\kk),  
 m^{\mu^{(i)}}_{{\lambda^{(i)}}',{\lambda^{(i)}}''} >0 
\}.
\end{equation}
For $\kk,\dd ,\m\in \Z_{>0}^\ell$, we denote
\begin{equation}
\label{eqn:set-of-irred}
\mathcal{I}(\kk,\dd,\mm) :=
\bigcup_{\pmb{\lambda}  \in \Par(\kk, (2\dd)^\mm)} \mathcal{I}(\pmb{\lambda}).
\end{equation}
If $\ell=1$, $\kk=(k),\dd= (d),\mm=(m)$, we will denote $\mathcal{I}(\kk,\dd,\mm)$ by
$\mathcal{I}(k,d,m)$.
Notice that for $\kk,\dd ,\m\in \Z_{>0} ^\ell$, 
\begin{eqnarray}
\label{eqn:set-of-irred2}
\mathcal{I}(\kk,\dd,\mm) &=& \prod_{i=1}^{\ell} \mathcal{I}(k_i,d_i,m_i).
\end{eqnarray}
\end{definition}

It follows directly from Part \eqref{item:prop:multiplicity1} of 
Proposition \ref{prop:multiplicity} that:

\begin{proposition}
\label{prop:restriction}
For $\kk,\dd ,\m\in \Z_{>0}^\ell$, and $\pmb{\mu} \in  \mathcal{I}(\kk,\dd,\m)$,
\[
\rank(\pmb{\mu}) \leq (2\dd)^\m.
\]
\end{proposition}

\begin{proof}
The proposition follows from Part (\ref{item:prop:multiplicity1}) of Proposition \ref{prop:multiplicity} and definition of $\mathcal{I}(\kk,\dd,\m)$ 
(cf. Notation \ref{def:set-of-irred}).
\end{proof}

\begin{remark}
\label{rem:restriction}
Note  that Proposition \ref{prop:restriction} implies that the Young diagram for each  $\mu \in \mathcal{I}(k,d,m)$ is contained in the union of $(2d)^m$ rows and$(2d)^m$ columns. This is shown in Figure
\ref{fig:young} for fixed $d,m$ and large $k$. 
The shaded area inside the $k \times k$ sized box contains all possible Young diagrams of
partitions of $k$. The darker part contains the partitions belonging to $\mathcal{I}(k,d,m)$. 

Since for every $d \leq k$, it is clear that 
\[
\card(\{\mu \in \Par(k)| \rank(\mu) = d\}) \leq 2 k^d,
\]
it follows immediately from Proposition \ref{prop:restriction} that for every fixed $d,m$, 
$\card(\mathcal{I}(k,d,m))$ is bounded by a polynomial in $k$.

\begin{figure}

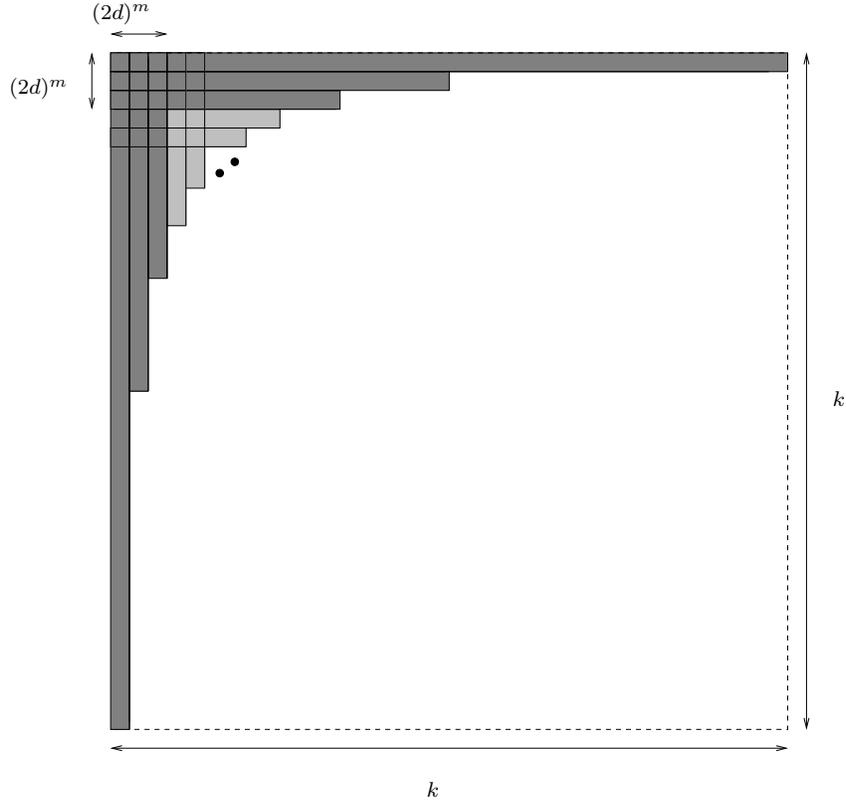

\caption{The shaded area contains all Young diagrams of partitions in $\Par(k)$, while the darker area contains the Young diagrams of the partitions in the subset $\mathcal{I}(k,d,m)\subset \Par(k)$ for fixed $d,m$ and large $k$.}
\label{fig:young}
\end{figure}
\end{remark}

The main structural result of this section is the following.

\begin{theorem}
\label{thm:main-product-of-symmetric}
Let $\kk=(k_1,\ldots,k_\ell),
\m =(m_1,\ldots,m_\ell), 
\dd=(d,\ldots,d) \in \Z_{>0}^\ell$, and $K = \sum_{i=1}^{\ell} m_i k_i$.
Let 
$P \in \R[\X^{(1)},\ldots,\X^{(\ell)}]$
be a 
non-negative
$\mathfrak{S}_\kk$-symmetric polynomial, 
with $\deg(P) \leq d$. Let  $V = \ZZ(P,\R^K)$.
Then, for all $\pmb{\mu} \in \Par(\kk)$,  $m_{\pmb{\mu}}(V,\F) > 0$ implies that
\begin{equation}
\label{eqn:restriction-on-specht}
\pmb{\mu} \in 
\mathcal{I}(\kk,\dd,\m).
\end{equation}
Moreover, for each  $\pmb{\mu} = (\mu^{(1)},\ldots,\mu^{(\ell)}) \in \mathcal{I}(\kk,\dd,\m)$,  
\begin{eqnarray}
\label{eqn:thm:main-product-of-symmetric}
  m_{\pmb{\mu}}(V,\F) &\leq& 
\sum_{\pmb{\lambda}=(\lambda^{(1)},\ldots,\lambda^{(\ell)})\in \Par(\kk, (2\dd)^\m)}
G(\pmb{\mu},\pmb{\lambda},\dd,\m),
\end{eqnarray}
where 
\begin{eqnarray}
\label{eqn:def-of-G}
G(\pmb{\mu},\pmb{\lambda},\dd,\m) &=&
\prod_{1 \leq i \leq \ell}
\left(
 (2 d)^{(m_i\length(\lambda^{(i)}))}
 \max_{\lambda^{(i)} = {\lambda^{(i)}}' \coprod {\lambda^{(i)}}''}
 m^{\mu^{(i)}}_{{\lambda^{(i)}}',{\lambda^{(i)}}''}
 \right)
 \end{eqnarray}
(the maximum on the right hand side is taken over all decompositions
$\lambda^{(i)} = {\lambda^{(i)}}' \coprod {\lambda^{(i)}}''$). 
\end{theorem}

\begin{proof}[Proof of Theorem \ref{thm:main-product-of-symmetric}]
We first assume that $V$ is bounded.
We replace $V$ by the set $S$ defined as the union of semi-algebraically connected components
of the set defined by 
$\Def(P,d',\zeta) \leq 0$
which are bounded over $\R$,
where $d'$ is the least even number such that $d' > d$ and
where $d'-1$ is prime. It follows from Bertrand's postulate that $d' \leq 2d$.
Using 
Proposition \ref{prop:alg-to-semialg},
Proposition \ref{prop:non-degenerate}, 
Lemma   \ref{lem:equivariant_morseA},
Proposition \ref{prop:equivariant-morseb}, 
and Proposition \ref{prop:half-degree}, 
we see that for $\mathbb{S}^{\pmb{\mu}}, \pmb{\mu} \in \Par(\kk)$ to occur with positive multiplicity in the isotypic decomposition of 
$\HH^*(V,\F)$, there must exist
$\pmb{\lambda}\in \Par(\kk, (2\dd)^{\mm})$
such that  $\pmb{\mu}  \in \overline{\Par}(\pmb{\lambda})$,
where $\overline{\Par}(\pmb{\lambda})$ is defined in Notation \ref{not:induced-rep-and-mult}. 
Eqn. \eqref{eqn:restriction-on-specht} now follows from 
Eqn. \eqref{eqn:set-of-irred},
and
inequality \eqref{eqn:thm:main-product-of-symmetric} follows from Part \eqref{item:prop:multiplicity3} of Proposition \ref{prop:multiplicity}.

More generally, introduce a new block of one variable, $Y$, and let 
\begin{eqnarray*}
\widetilde{P} &=&  P  +  \left(\eps^2 \cdot \sum_{1\leq h \leq \ell} \sum_{1 \leq i \leq m_h} \sum_{1 \leq j \leq k_h} (X_{i,j}^{(h)})^2 + Y^2 - 1\right)^2,\\
\widetilde{Q} &=&  P  +  \left(\eps^2 \cdot \sum_{1\leq h \leq \ell} \sum_{1 \leq i \leq m_h} \sum_{1 \leq j \leq k_h} (X_{i,j}^{(h)})^2 - 1\right)^2,
\end{eqnarray*}
and let 
\begin{eqnarray*}
\widetilde{V} &=&  \ZZ(\widetilde{P},\R\la\eps\ra^{K+1}), \\
\widetilde{W} &=&  \ZZ(\widetilde{Q},\R\la\eps\ra^{K}), \\
\widetilde{T} &=&  \Ext(V,\R\la\eps\ra) \cap \overline{B_K(0,1/\eps)}.
\end{eqnarray*}
We let $\kk' = (\kk,1)$. Clearly, $\widetilde{V}$ is $\mathfrak{S}_{\kk'}$-symmetric and bounded over $\R\la\eps\ra$, and 
$\widetilde{W}$ is $\mathfrak{S}_{\kk}$-symmetric and also bounded over $\R\la\eps\ra$.

It follows from the conical structure theorem at infinity for semi-algebraic sets, that:
\begin{enumerate}[i]
\item
$\Ext(V,\R\la\eps\ra)$ is semi-algebraically homeomorphic to $\widetilde{T}$; and
\item
$\widetilde{V} = \widetilde{V}_{+}\cup \widetilde{V}_{-}$,  where for $\sigma \in \{+,-\}$,
$\widetilde{V}_{\sigma}$ is the intersection of $\widetilde{V}$ with the half-space defined by $\sigma Y \geq 0$;
\item
$\widetilde{T}$, and hence $\Ext(V,\R\la\eps\ra)$,  is semi-algebraically homeomorphic to each of $\widetilde{V}_{+}, \widetilde{V}_{-}$;
\item
$\widetilde{W} = \widetilde{V}_{+}\cap \widetilde{V}_{-}$.
\end{enumerate}

It follows from Proposition \ref{prop:MV} that 
\[
m_{\pmb{\mu}}(V,\F) \leq \frac{1}{2} (m_{\pmb{\mu}'}(\widetilde{V},\F) + m_{\pmb{\mu}}(\widetilde{W},\F)).
\]
The theorem now follows from the bounded case proved before, noticing that the result in the bounded case implies that,
\[
m_{\pmb{\mu}'}(\widetilde{V},\F), m_{\pmb{\mu}}(\widetilde{W},\F)
\] 
are both bounded by 
\begin{eqnarray*}
\sum_{\pmb{\lambda}=(\lambda^{(1)},\ldots,\lambda^{(\ell)})\in \Par(\kk, (2\dd)^\m)}
G(\pmb{\mu},\pmb{\lambda},\dd,\m).
\end{eqnarray*}
\end{proof}

\subsection{Proofs of Theorems \ref{thm:main-product-of-symmetric-quantitative}, \ref{thm:equivariant}, and \ref{thm:main-product-of-symmetric-quantitative-complex}}
\begin{proof}[Proof of Theorem \ref{thm:main-product-of-symmetric-quantitative}]
Theorem \ref{thm:main-product-of-symmetric-quantitative} follows from
Theorem \ref{thm:main-product-of-symmetric}, Proposition \ref{prop:restriction}, and 
Part \eqref{item:prop:multiplicity3}
of
Proposition \ref{prop:multiplicity}.
\end{proof}

\begin{proof}[Proof of Theorem \ref{thm:equivariant}]
Theorem \ref{thm:equivariant} follows from
Theorem \ref{thm:main-product-of-symmetric} and 
Remark \ref{rem:multiplicity}.
\end{proof}

\begin{proof}[Proof of Theorem \ref{thm:main-product-of-symmetric-quantitative-complex}]
Substituting  $\X^{(j)} = \Y^{(j)}+ i \ZB^{(j)}, 1 \leq j \leq \ell$ in $\mathcal{P}$
and separating the real and imaginary parts, obtain another family of polynomials,
$\mathcal{Q} \subset \R[\Y^{(1)},\ZB^{(1)},\ldots, \Y^{(\ell)},\ZB^{(\ell)}]$ with
$\deg_{\Y^{(j)}}(Q),\deg_{\ZB^{(j)}}(Q) \leq d, 1 \leq j \leq \ell$, such that the polynomials in
$\mathcal{Q}$ are $\mathfrak{S}_{\mathbf{k}}$-symmetric.

Now apply Theorem \ref{thm:main-product-of-symmetric-quantitative} with
$\kk=(k_1,\ldots,k_\ell)$,
$\m =(2m_1,\ldots,2m_\ell)$, and  
$\dd=(d,\ldots,d)$.
\end{proof}

\subsection{Proof of Theorem \ref{thm:main-product-of-symmetric-sa-quantitative}}

We first prove a more structural result from which Theorem \ref{thm:main-product-of-symmetric-sa-quantitative} will follow easily.

\begin{theorem}
\label{thm:main-product-of-symmetric-sa}
Let $\kk=(k_1,\ldots,k_\ell),
\m =(m_1,\ldots,m_\ell), 
\dd=(d,\ldots,d) \in \Z_{>0}^\ell$, and $K = \sum_{i=1}^{\ell} k_i m_i$.
Let 
$\mathcal{P} \subset \R[\X^{(1)},\ldots,\X^{(\ell)}]^{\mathfrak{S}_\kk}_{\leq \dd}$
be a finite set of 
of polynomials, and let $\card(\mathcal{P}) = s$. Let  
$S \subset \R^K$  be a $\mathcal{P}$-closed semi-algebraic set.

Then, for all $\pmb{\mu} \in \Par(\kk)$,  $m_{\pmb{\mu}}(S,\F) > 0$ implies that
\[
\pmb{\mu} \in 
\mathcal{I}(\kk,\dd,\m).
\]
Moreover, 
let $D=D (\mathbf{k},\mathbf{m},d) = \sum_{i=1}^{\ell} \min
  (m_i k_i ,  d^{m_i})$.
Then,
for each  
\[
\pmb{\mu} = (\mu^{(1)},\ldots,\mu^{(\ell)}) \in \mathcal{I}(\kk,\dd,\m),
\]  
\begin{eqnarray*}
  m_{\pmb{\mu}}(S,\F)  &\leq& 
    \sum_{i=0}^{D-1}
    \sum_{j=1}^{D-i} \binom{2 s+1}{j} 6^{j} \cdot \left(\sum_{\pmb{\lambda}=(\lambda^{(1)},\ldots,\lambda^{(\ell)})\in \Par(\kk, (4\dd)^\m)}
G(\pmb{\mu},\pmb{\lambda},2 \dd,\m)\right),
\end{eqnarray*}
where 
$G(\pmb{\mu},\pmb{\lambda},\dd,\m)$ is defined in Eqn. \eqref{eqn:def-of-G}.
\end{theorem}

Before proving Theorem \ref{thm:main-product-of-symmetric-sa} 
we first need a few preliminary definitions and results.

\begin{definition}[$\ell$-general position]
  For any finite family $\mathcal{P} \subset \R [ X_{1} , \ldots ,X_{k} ]$ and
  $\ell \geq 0$, we say that $\mathcal{P}$ is in $\ell$-general position with
  respect to a semi-algebraic set $V \subset \R^{k}$ if for any subset
  $\mathcal{P}' \subset \mathcal{P}$, with $\card (\mathcal{P}') > \ell$, $\ZZ (\mathcal{P}' ,V) = \emptyset$. 
\end{definition}

Let $\kk=(k_1,\ldots,k_\ell),
\m =(m_1,\ldots,m_\ell)\in  \Z_{>0}^\ell$, and $K = \sum_{i=1}^{\ell} k_i m_i$.
Let 
$\mathcal{P} = \{P_1,\ldots, P_s\} \subset \R[\X^{(1)},\ldots,\X^{(\ell)}]^{\mathfrak{S}_\kk}_{\leq d}$
be a finite set of polynomials, and
let  
$S \subset \R^K$  be a $\mathcal{P}$-closed semi-algebraic set.
Let $\overline{\eps} = \left(\eps_{1} , \ldots , \eps_{s} \right)$
be a tuple of new variables, and let $\mathcal{P}_{\overline{\eps}} =
\bigcup_{1 \leq i \leq s} \left\{ P_{i}   \pm \eps_{i} \right\}$. We have the
following two lemmas.

\begin{lemma}
  \label{lem:gen-pos1-with-parameters}Let
  \begin{eqnarray*}
    D (\mathbf{k},\mathbf{m},d) & = & \sum_{i=1}^{\ell} \min (k_{i}m_i ,d^{m_i}) .
  \end{eqnarray*}
  The family $\mathcal{P}_{\overline{\eps}} \subset \R'[\X^{(1)} , \ldots ,\X^{(\ell)}]$ is in $D$-general position with respect to any semi-algebraic subset $Z' \subset \R'^K$, 
where $\R' = \R \langle \overline{\eps} \rangle$ (cf. Notation \ref{not:Puiseux}), and
where $Z' = \Ext(Z,\R'^K)$ (cf. Notation \ref{not:Ext}), and $Z\subset \R^{K}$ is a semi-algebraic set stable under the action of $\mathfrak{S}_{\mathbf{k}}$.
  \end{lemma}

\begin{proof}
The lemma follows from the fact that the ring of multi-symmetric polynomials is generated by the
multi-symmetric power sum polynomials \cite[Theorem 1.2]{Dalbec}, and the cardinality of the
set of multi-symmetric power sum polynomials in the variables $X^{(i)}$ of degree bounded by $d$
is bounded by $d^{m_i}$.
\end{proof}

Let $\Phi$ be a $\mathcal{P}$-closed formula, and let $S= \RR (\Phi ,V)$ be
bounded over $\R$.
\begin{notation}[Multiplicities]
For $\pmb{\mu} \in \Par(\kk)$ and $i \geq 0$, we will denote
\begin{eqnarray*}
m_{i,\pmb{\mu}}(\Phi,\F) &=&  m_{i,\pmb{\mu}}(S,\F), \\
m_{\pmb{\mu}}(\Phi,\F) &=&  m_{\pmb{\mu}}(S,\F).
\end{eqnarray*}
\end{notation}
  
Let $\Phi_{\overline{\eps}}$ be the
$\mathcal{P}_{\overline{\eps}}$-closed formula obtained from $\Phi$ be
replacing for each $i,1 \leq i \leq s$,
\begin{enumerate}[i.]
  \item each occurrence of $P_{i} \leq 0$ by $P_{i} - \eps_{i} \leq 0
  $, and
  
  \item each occurrence of $P_{i} \geq 0$ by $P_{i} + \eps_{i} \geq 0
  $.
\end{enumerate}
Let $\R' = \R \left\langle \eps_{1} , \ldots , \eps_{s} \right\rangle$, and 
$S_{\overline{\eps}} = \RR(\Phi_{\overline{\eps}} , \R'^{K})$.

\begin{lemma}
  \label{lem:gen-pos2-with-parameters} For any $r>0$, $r \in \R$, the
  semi-algebraic set set $\Ext (S \cap \overline{B_{K} (0,r)} , \R')$ is contained in 
  $S_{\overline{\eps}} \cap \overline{B_{K} (0,r)}$, and the inclusion 
  $\Ext(S \cap \overline{B_{K} (0,r)} , \R') \hookrightarrow S_{\overline{\eps}} \cap
  \overline{B_{K} (0,r)}$ is a semi-algebraic homotopy equivalence.
  The induced isomorphism,
  \[
  \HH(S_{\overline{\eps}} \cap \overline{B_{K} (0,r)},\F) \overset{\sim}{\rightarrow} \HH^*(\Ext(S \cap \overline{B_{K} (0,r)} , \R'),\F) 
  \]
  is an isomorphism of $\mathfrak{S}_{\mathbf{k}}$-modules.
\end{lemma}

\begin{proof} The proof is similar to the one of Lemma 16.17 in
{\cite{BPRbook2}}.
\end{proof}

\begin{remark}
  \label{rem:gen-pos3} In view of Lemmas \ref{lem:gen-pos1-with-parameters} and
  \ref{lem:gen-pos2-with-parameters} we can assume (at the cost of doubling the number of
  polynomials) after possibly replacing $\mathcal{P}$ by
  $\mathcal{P}_{\overline{\eps}}$, and $\R$ by $\R \left\langle \eps_{1} ,
  \ldots , \eps_{s} \right\rangle$, that the family $\mathcal{P}$ is in
  $D(\mathbf{k},\mathbf{m},d)$-general position.
\end{remark}

Now, let $\delta_{1} , \cdots , \delta_{s}$ be new infinitesimals, and let
$\R' = \R \langle \delta_{1} , \ldots , \delta_{s} \rangle$.

\begin{notation}[Infinitesimal thickening]
  We define $\mathcal{P}_{>i} = \{P_{i+1} , \ldots ,P_{s} \}$ and
  \begin{eqnarray*}
    \Sigma_{i} & = & \{P_{i} =0,P_{i} = \delta_{i} ,P_{i} = - \delta_{i}
    ,P_{i} \geq 2 \delta_{i} ,P_{i} \leq -2 \delta_{i} \} ,\\
    \Sigma_{\le i} & = & \{\Psi \mid \Psi = \bigwedge_{j=1, \ldots ,i}
    \Psi_{i} , \Psi_{i} \in \Sigma_{i} \} .
  \end{eqnarray*}
  Note that for each $\Psi \in \Sigma_{i}$, $\RR(\Psi , \R \langle
  \delta_{1} , \ldots , \delta_{i} \rangle^{K})$ is symmetric with respect
  to the action of $\mathfrak{S}_{\kk}$,  and for
  $\Psi \neq \Psi'$, $\Psi , \Psi' \in \Sigma_{\leq i}$,
  \begin{eqnarray}
    \RR \left(\Psi , \R \langle \delta_{1} , \ldots , \delta_{i} \ra^{K} \right)
    \cap \RR \left(\Psi' , \R \la \delta_{1} , \ldots , \delta_{i} \ra^{K}
    \right) & = & \emptyset .  \label{eqn:disjoint}
  \end{eqnarray}

  If $\Phi$ is a $\mathcal{P}$-closed formula, we denote
  \begin{eqnarray*}
    \RR_{i} (\Phi) & = & \RR(\Phi , \R \la \delta_{1} , \ldots ,\delta_{i} \ra^{K}) ,
  \end{eqnarray*}
  and
  \begin{eqnarray*}
    \RR_{i} (\Phi \wedge \Psi) & = & \RR(\Psi , \R \la \delta_{1} ,
    \ldots , \delta_{i} \ra^{K}) \cap \RR_{i} (\Phi) .
  \end{eqnarray*}
\end{notation}

The proof of the following proposition is very similar to Proposition 7.39 in
{\cite{BPRbook2}} where it is proved in the non-symmetric case.

\begin{proposition}
  \label{7:prop:closed-with-parameters}
  For every $\mathcal{P}$-closed formula
  $\Phi$, and $\pmb{\mu} \in \Par(\kk)$,
   such that $\RR(\Phi)$ is bounded,
  \[  m_{\pmb{\mu}}(\Phi,\F) \leq
     \sum_{\substack{
       \Psi \in \Sigma_{\le s}\\
       \RR_{s} (\Psi , \R'^{K}) \subset \RR_{s} (\Phi , \R'^{K})
       }}
      m_{\pmb{\mu}}(\Psi,\F) . 
     \]
\end{proposition}

\begin{proof}
 The symmetric spaces $\RR \left(\Psi , \Ext \left(V, \R' \right) \right) ,
\Psi \in \Sigma_{\leq s}$ are disjoint by (\ref{eqn:disjoint}). The 
proposition now follows from Schur's lemma,
and the proof of Proposition 7.39 in {\cite{BPRbook2}}.
\end{proof}

\begin{proposition}
  \label{7:prop:betti closed}
  Suppose for $\pmb{\mu} \in \Par(\kk)$ and $i \geq 0$,
  $m_{i, \pmb{\mu}}(S,\F) >0$. 
 Then, 
 \begin{equation}
\label{eqn:restriction-on-specht-sa}
\pmb{\mu}  \in 
\mathcal{I}(\kk,\dd,\m),
\end{equation}
where $\dd = (d,\ldots,d)$.
For $i \geq 0$, and $\pmb{\mu}  \in 
\mathcal{I}(\kk,\dd,\m)$,
  \[ 
  \sum_{\Psi \in \Sigma_{\le s}} m_{i,\pmb{\mu}}(\Psi,\F) \leq 
  \sum_{j=0}^{D(\kk,\mm,d)} \binom{s}{j} 6^{j} 
  F(\pmb{\mu}, \kk,\mm, 2d),
  \]
  where 
  \begin{equation}
  \label{eqn:definition-of-F}
  F(\pmb{\mu},\kk,\mm,d) = 
  \sum_{\pmb{\lambda}=(\lambda^{(1)},\ldots,\lambda^{(\ell)})\in \Par(\kk, (2\dd)^\m)}
G(\pmb{\mu},\pmb{\lambda},\dd,\m),
 \end{equation}
 and 
$G(\pmb{\mu},\pmb{\lambda},\dd,\m)$ is defined in \eqref{eqn:def-of-G}.
 \end{proposition}

In order to prove Proposition \ref{7:prop:betti closed} we first need the following lemmas.

Let for $1 \leq i \leq s$, 
$Q_{i} =P_{i}^{2} (P_{i}^{2} - \delta_{i}^{2})^{2} (P_{i}^{2} -4
\delta_{i}^{2})$.

For $j \ge 1$ let,
\begin{eqnarray*}
  V'_{j} & = & \RR(\bigvee_{1 \leq i \leq j} Q_{i} =0, \R \la\delta_{1} , \ldots , \delta_{j} \ra^{K}) ,\\
  W'_{j} & = & \RR(\bigvee_{1 \leq i \leq j} Q_{i} \geq 0, \R \la \delta_{1} , \ldots , \delta_{j} \ra^{K}) .
\end{eqnarray*}
\begin{lemma}
  \label{lem:gen-position} Let $I \subset [ 1,s ]$, $\sigma  =  (\sigma_{1} ,
  \ldots , \sigma_{s}) \in \{ 0, \pm 1, \pm 2 \}^{s}$ and let
  $\mathcal{P}_{I, \sigma} =  \bigcup_{i \in I} \{ P_{i} + \sigma_{i}
  \delta_{i} \}$. Then, $\ZZ \left(P_{I, \sigma} , \R'^K \right) = \emptyset$,
  whenever $\card(I) >D$.
\end{lemma}

\begin{proof}
This follows from the fact that $\mathcal{P}$ is in
$D$-general position by Remark
\ref{rem:gen-pos3}.
\end{proof}

\begin{lemma}
  \label{7:lem:union2} For each $\pmb{\mu}  \in 
\mathcal{I}(\kk,\dd,\m)$, and $i \geq 0$,
\[ m_{i,\pmb{\mu}} (V'_{j},\F) \leq 
 \sum_{p=1}^{\min (j,D)} \binom{j}{p} 5^{p} 
F (\pmb{\mu}, \mathbf{k},\m,2d) 
\]
(see  \eqref{eqn:definition-of-F} above for the definition of 
$F (\pmb{\mu}, \mathbf{k},\m,d)$).
\end{lemma}

\begin{proof}

The set $\RR ((P_{j}^{2} (P_{j}^{2} - \delta_{j}^{2}
)^{2} (P_{j}^{2} -4 \delta_{j}^{2})=0), \R \la \delta_{1} , \ldots ,
\delta_{j} \ra^{K})$ is the disjoint union of

\begin{equation}
 \begin{array}{c}
     \RR (P_{i} =0, \R \la \delta_{1} , \ldots , \delta_{j} \ra^{K}) ,\\
     \RR (P_{i} = \delta_{i} , \R \la \delta_{1} , \ldots , \delta_{j} \ra^{K}
    ) ,\\
     \RR (P_{i} =  -\delta_{i} , \R \la \delta_{1} , \ldots , \delta_{j}
     \ra^{K}) ,\\
     \RR (P_{i} =2 \delta_{i} , \R \la \delta_{1} , \ldots , \delta_{j}
     \ra^{K}) ,\\
     \RR (P_{i} =  -2 \delta_{i} , \R \la \delta_{1} , \ldots , \delta_{j}
     \ra^{K}) .
   \end{array} \label{eqn:list}
\end{equation}

It follows from part (\ref{7:prop:prop1:itema}) of Proposition \ref{7:prop:prop1}  
that $m_{i,\pmb{\mu}} (V'_{j} ,\F)$ 
is bounded by the sum for 
$1 \leq p \leq i+1$, of the multiplicities of  $\mathbb{S}^{\pmb{\mu}}$
in the $(i- p +1)$-th cohomology module of all possible
non-empty sets obtained by the intersection of $p$ distinct sets from amongst
amongst the sets listed in \eqref{eqn:list}.
Because of the fact that the set of polynomials $\mathcal{P}$ is in $D$-general position it follows that
all such intersections will be empty if $p > D$ or $p > j$. Moreover, Thus, the total number of non-empty
intersections that we need to consider is bounded by 
\[
\sum_{p=1}^{\min(j,D)} \binom{j}{p} 5^p.
\] 

It now follows from 
Theorem 
\ref{thm:main-product-of-symmetric}
applied to the non-negative
symmetric polynomials 
$P_{i}^{2} , (P_{i} \pm\delta_{i})^{2} , (P_{i} \pm 2 \delta_{i})^{2}$, 
and noting that the
degrees of these polynomials are bounded by 
$2d$, 
that
\begin{eqnarray*}
  m_{i, \pmb{\mu}} (V'_{j},\F) 
  & \leq & 
  \sum_{p=1}^{\min (j,D)} \binom{j}{p} 5^{p} 
 F (\pmb{\mu}, \mathbf{k},\m,2d).
\end{eqnarray*}
\end{proof}

\begin{lemma}
  \label{7:lem:union1}
  For each $\pmb{\mu}  \in 
\mathcal{I}(\kk,\dd,\m)$, and $i \geq 0$,
  \[ m_{i,\pmb{\mu}} (W'_{j},\F) \leq \sum_{p=1}^{\min (
     j,D)} \binom{j}{p} 5^{p} 
     F (\pmb{\mu}, \mathbf{k},\m,2 d) +m_{i,\pmb{\mu}} (\R \la \delta_{1} , \ldots , \delta_{j} \ra^{K},\F) . \]
 \end{lemma}

\begin{proof} Let
\[ T = \RR \left(\bigwedge_{1 \leq i \leq j} Q_{i} \leq 0 \vee \bigvee_{1 \leq
   i \leq j} Q_{i} =0, \Ext (Z, \R \la \delta_{1} , \ldots , \delta_{i} \rangle)
   \right) . \]

Now, from the fact that 
\[ W'_{j} \cup T= \R \la \delta_{1} , \ldots ,
   \delta_{j} \ra^{k} ,W'_{j} \cap T=V'_{j},\] 
and Proposition \ref{prop:MV} it follows that
\begin{eqnarray*}
  m_{i,\pmb{\mu}} (W'_{j} ,\F) & \leq & m_{i,\pmb{\mu}} ((
  W'_{j} \cap T),\F) +m_{i,\pmb{\mu}} ((W'_{j}
  \cup T),\F)\\
  & = & m_{i,\pmb{\mu}} (V'_{j},\F) + m_{i,\pmb{\mu}} (
  \R \la \delta_{1} , \ldots , \delta_{j} \ra^{K},\F).
\end{eqnarray*}
We conclude using Lemma \ref{7:lem:union2}.
\end{proof}

\begin{proof}[Proof of Proposition \ref{7:prop:betti closed}]
Using part (\ref{7:prop:prop1:itemb}) of Proposition \ref{7:prop:prop1} we get that
\begin{eqnarray*}
  \sum_{\Psi \in \Sigma_{\le s}} m_{i,\pmb{\mu}} (\Psi,\F) & \leq & \sum_{j=1}^{\min(D,K-i)}
  \sum_{\begin{array}{c}
    J \subset \{ 1, \ldots ,s \}\\
    \card(J) =j
  \end{array}} m_{i+j-1,\pmb{\mu}} (S^{J},\F) \\
  && + \binom{s}{K-i}
  m_{K,\pmb{\mu}} (S^{\emptyset},\F).
\end{eqnarray*}
It follows from Lemma \ref{7:lem:union1} that,
\begin{eqnarray*}
  m_{i+j-1,\pmb{\mu}} (S^{J}) & \leq & \sum_{p =1}^{\min (j,D)} \binom{j}{p}
  5^{p} 
 F (\pmb{\mu},\mathbf{k},\m,2 d) +m_{K,\pmb{\mu}}(\R^{K},\F).
\end{eqnarray*}
Hence,
\begin{eqnarray*}
  \sum_{\Psi \in \Sigma_{\le s}} m_{i,\pmb{\mu}} (\Psi,\F)
  & \leq & \sum_{j=1}^{D}
  \sum_{\substack{
    J \subset \{ 1, \ldots ,s \}\\
    \card (J) =j}}
 m_{i+j-1,\pmb{\mu}} (S^{J}, \F) + \binom{s}{K-i}
  m_{K,\pmb{\mu}} (S^{\emptyset},\F)\\
  & \leq & \sum_{j=1}^{D} \binom{s}{j} \left(\sum_{p=1}^{\min (
  j,D)} \binom{j}{p} 5^{p}  
  F (\pmb{\mu},\mathbf{k},\m,2d) \right)\\
  & \leq & \sum_{j=1}^{D} \binom{s}{j} 6^{j} 
  F (\pmb{\mu},\mathbf{k},\m,2d).
\end{eqnarray*}
\end{proof}

\begin{proof}[Proof of Theorem \ref{thm:main-product-of-symmetric-sa}]
We first add an extra polynomial, $\delta(X_1^2+\cdots+X_K^2) -1$ to the set $\mathcal{P}$, replace the field $\R$, by
$\R\la\delta\ra$, and replace the given formula $\mathcal{P}$-closed formula $\Phi$ by the formula
$\Phi \wedge ( \delta(X_1^2+\cdots+X_K^2) -1\leq 0)$. Notice that the new set ${\rm Reali}(\Phi)$
is bounded in $\R\la\delta\ra^k$ and is $\mathfrak{S}_{\kk}$-symmetric.
The theorem now follows from
Propositions \ref{7:prop:closed-with-parameters} and \ref{7:prop:betti closed}.
\end{proof}

\begin{proof}[Proof of Theorem \ref{thm:main-product-of-symmetric-sa-quantitative}]
Follows immediately from Theorem \ref{thm:main-product-of-symmetric-sa} and 
Proposition \ref{prop:multiplicity}. 
\end{proof}

\subsection{Proof of Theorem \ref{thm:symmetric-complex-projective}}
\begin{proof}[Proof of Theorem \ref{thm:symmetric-complex-projective}]
Let $\Sphere^{2k+1} \subset \C^{k+1}$ denote the unite sphere defined by 
$|Z_0|^2 +\cdots + |Z_k|^2 =1$. 
Consider the Hopf fibration $\phi: \Sphere^{2k+1} \rightarrow \PP_\C^k$, defined by
$(z_0,\ldots,z_k) \mapsto (z_0:\cdots:z_k)$. 
We denote by $\tilde{V} = \phi^{-1}(V)$.
We have the following commutative diagram:
\[
\xymatrix{
\tilde{V} \ar[r]^{i} \ar[d]^{\phi|_{\tilde{V}}} & \Sphere^{2k+1} \ar[d]^\phi \\
V \ar[r]^i & \PP_\C^k
}
\]
Note that $\tilde{V}$ is a $\Sphere^1$-bundle over $V$, and using the fact that $\PP_\C^k$ is simply connected, 
there is a $\mathfrak{S}_{k+1}$-equivariant spectral sequence degenerating at its $E_3$ term converging to the cohomology of  $\tilde{V}$. 

The $E_3$-term of the spectral sequence is given by
\begin{eqnarray*}
E_2^{p,q} & \cong & \HH_p(V,\F), \mbox{ if } q=0,1,\\
E_2^{p,q} &=& 0, \mbox{ else },
\end{eqnarray*} 
and the differentials $d_2^{p,q}: E_2^{p,q} \rightarrow E_2^{p+2,q-1}$ shown below. 

{\tiny
\[
\xymatrix{
 \vdots & \vdots & \vdots & \vdots & \vdots & \vdots & \vdots  \\
 0 \ar[rrd]^{d_2^{-1,2}} &0 &0 &\cdots &0\ar[rrd]^{d_2^{i,2}}  &0 & 0\\
 0 \ar[rrd]^{d_2^{-1,1}} &\HH^0(V,\F) & \HH^1(V,\F)   &\cdots& \HH^i(V,\F)\ar[rrd]^{d_2^{i,1}}  & \HH^{i+1}(V,\F) & \HH^{i+2}(V,\F)\\
0  & \HH^0(V,\F) & \HH^1(V,\F)  &\cdots
& \HH^i(V,\F)  & \HH^{i+1}(V,\F) &  \HH^{i+2}(V,\F) \\
}
\]
}
Fix $\lambda \vdash k+1$, and recall that we denote for each $i \geq 0$, $m_{i,\lambda}(V,\F)$ (resp. ${m}_{i,\lambda}(\tilde{V},\F)$) the multiplicity of
$\mathbb{S}^\lambda$ in $\HH^i(V,\F)$ (resp. $\HH^i(\tilde{V},\F)$).

We observe that since $\HH^0(V,\F) \cong_{\mathfrak{S}_{k+1}} \HH^0(\tilde{V},F)$, 
we have for all $\lambda \vdash k+1$,
\begin{equation}
\label{eqn:spectral-multiplicity-0}
m_{0,\lambda}(V,\F) = {m}_{\lambda,0}(\tilde{V},\F).
\end{equation}

Also, note that it follows from the fact that the spectral sequence $E_r^{p,q}$  degenerates at its
$E_3$ term that, 
\[
\HH^1(V,\F) \oplus \Ker(d_2^{0,1}) \cong_{\mathfrak{S}_{k+1}}  \HH^1(\tilde{V},\F),
\]
and we obtain from the fact that the spectral sequence $E^r_{p,q}$ is $\mathfrak{S}_{k+1}$-equivariant that
\begin{equation}
\label{eqn:spectral-multiplicity-1}
m_{1,\lambda}(V,\F) \leq {m}_{1,\lambda}(\tilde{V},\F).
\end{equation}

More generally, we have from the $E^2$-term of the spectral sequence that
\[
\HH^i(\tilde{V},\F)  \cong_{\mathfrak{S}_{k+1}}  \mathrm{coker}(d_2^{i-2,1}) \oplus \Ker(d_2^{i-1,1}).
\]

For $\lambda \vdash k+1, i\geq 0$, 
and any finite dimensional $\F$-representation $W$ of $\mathfrak{S}_{k+1}$, we denote 
by $\mult_{\lambda}(W,\F)$ the multiplicity of $\mathbb{S}^\lambda$ in $W$.

Since, 
\[
\HH^i(V,\F) \cong_{\mathfrak{S}_{k+1}}  \mathrm{Im}(d_2^{i-2,1}) \oplus \mathrm{coker}(d_2^{i-2,1}),
\]
we have for all $\lambda \vdash k+1, i \geq 0$,
\[
\mult_{\lambda}( \mathrm{coker}(d_2^{i-2,1}),\F) = m_{i,\lambda}(V,\F)  -  \mult_{\lambda}(\mathrm{Im}(d_2^{i-2,1}),\F),
\]
and we also have for $i \geq 2$,

\begin{equation}
\label{eqn:Im}
\mult_\lambda(\mathrm{Im}(d_2^{i-2,1}),\F) \leq m_{i-2,\lambda}(V,\F).
\end{equation}
This implies that for all $\lambda \vdash k+1, i \geq 0$
\[
m_{i,\lambda}(\tilde{V},\F) =  (m_{i,\lambda}(V,\F) - \mult_{\lambda}(\mathrm{Im}(d_2^{i-2,1}),\F)) + 
\mult_\lambda(\mathrm{coker}(d_2^{i-1,1}),\F).
\]
It follows that
\begin{eqnarray*}
m_{i,\lambda}(V,\F) &=& 
m_{i,\lambda}(\tilde{V},\F) + \mult_{\lambda}(\mathrm{Im}(d_2^{i-2,1}),\F)) - 
\mult_\lambda(\mathrm{coker}(d_2^{i-1,1}),\F) \\
                     &\leq& m_{i,\lambda,i}(\tilde{V},\F) + \mult_\lambda(\mathrm{Im}(d_2^{i-2,1}),\F) \\
                     &\leq& m_{i,\lambda}(\tilde{V},\F) + m_{i-2,\lambda}(V,\F) \mbox{ using \eqref{eqn:Im}}.
                     \end{eqnarray*}

Finally we have shown that for each $\lambda \vdash k+1$ and $i \geq 2$,
 
\begin{eqnarray}
\label{eqn:spectral-multiplicity-general}
 m_{i,\lambda}(V,\F) &\leq& m_{i,\lambda}(\tilde{V},\F) + m_{i-2,\lambda}(V,\F)  \nonumber\\
 				&\leq& \sum_{0\leq j \leq \lfloor \frac{i}{2}\rfloor } m_{i - 2j,\lambda}(\tilde{V},\F) \mbox{ using induction}. 
\end{eqnarray}
 
The theorem follows from applying Theorem \ref{thm:main-product-of-symmetric-quantitative} to the set $\tilde{V}$, and 
inequalities \eqref{eqn:spectral-multiplicity-0}, \eqref{eqn:spectral-multiplicity-1}, and
\eqref{eqn:spectral-multiplicity-general}.
\end{proof}

\begin{remark}
Note that in the proof of Theorem \ref{thm:symmetric-complex-projective} it is possible to replace the spectral sequence argument altogether
by an argument using the equivariant version of the Gysin exact sequence. 
\end{remark}

\subsection{Proof of Theorem \ref{thm:descent2-quantitative-new}}
\begin{proof}[Proof of Theorem \ref{thm:descent2-quantitative-new}]
Let $P^{(p)} \in \R[\X,\Y_0,\ldots,\Y_{p}]$ be defined by
\[
 P^{(p)} = P(\X,\Y_0)+ \cdots+ P(\X,\Y_{p}),
 \]
\[
V^{(p)} = \ZZ(P^{(p)},\R^{k+(p+1)m}),
\]
 and 
\[
\kk(p) = (\underbrace{1,\ldots,1}_{k},p+1).
\]

  Notice that since $V$ is bounded, so is $V^{(p)}= \ZZ(P,\R^{k+(p+1)m})$, and moreover,
 $V^{(p)}$ is semi-algebraically homeomorphic to $W^{(p)}_\pi(V)$ (cf. Eqn. \eqref{eqn:definition-of-fibered-product}).
 Moreover, 
 $\deg(P^{(p)}) = \deg(P)$, and $P^{(p)}$ is symmetric in $(\Y_0,\ldots,\Y_{p})$, and is thus
 $\mathfrak{S}_{\kk(p)}$-symmetric.
 
By Theorem \ref{thm:descent2},
\begin{eqnarray*}
    b (\pi (V) ,\F) & \leq & \sum_{0 \leq p<k}
    b_{\mathfrak{S}_{\kk(p)}}(V^{(p)},\F).
\end{eqnarray*}

Now using Theorem  \ref{thm:equivariant},
\[
b_{\mathfrak{S}_{\kk(p)}}(V^{(p)},\F)
\leq (p+1)^{(2d)^{m}} (O(d))^{k+m (2d)^m + 1},
\]
and hence,
\begin{eqnarray*}
  b (\pi (V) ,\F) & \leq & \sum_{0 \leq p<k} 
  b_{\mathfrak{S}_{\kk(p)}}(V^{(p)},\F) \\
  &\leq & \sum_{0 \leq p<k} (p+1)^{(2d)^{m}} (O(d))^{k+ m (2d)^m +1} \\
  &\leq & k^{(2d)^{m}} (O(d))^{k+m (2d)^m +1}.
\end{eqnarray*}
 
This completes the proof of the theorem.
\end{proof}

\section{Conclusion and open problems}
\label{sec:conclusion}

In this paper we have proved polynomial bounds on the number and the multiplicities of the
irreducible representations of the symmetric group (or more generally product of symmetric groups) that appear in the cohomology modules of symmetric real algebraic and more generally real semi-algebraic sets. We have given several applications of the main results, including to improve existing
bounds on the topological complexity of sets defined as images of semi-algebraic maps, 
and proving lower bounds on the degrees etc. We end with some open problems and future research
directions.
 
\subsection{Representational Stability Question}
The bounds on the multiplicities that we prove in this paper are all polynomial in the number of variables
(for fixed  degrees).  Motivated by the recently developed theory of FI-modules \cite{Church-et-al} it makes sense  to ask whether it is possible to prove some stability result as $k \rightarrow \infty$.
We formulate one such question below.

Let $\kbb$ be a field, and let $A(\kbb)$ denote the polynomial ring $\kbb[(X_i)_{i \in \mathbb{N}} ]$ in the denumerable set of variables 
$\{X_1,X_2,\ldots\}$.

Let $\mathfrak{S}_\infty$ denote the infinite symmetric group, whose elements are bijections
$\mathbb{N} \rightarrow \mathbb{N}$ which keep all but finitely many elements of $\mathbb{N}$ fixed.

We say that an ideal  $I\subset A(\kbb)$ is symmetric it is stable under the natural action of $\mathfrak{S}_\infty$ 
permuting the variables. We say that a symmetric ideal $I \subset A(\kbb)$ is finitely generated, if there exists a finite
subset $\mathcal{F} \subset A(\kbb)$ such that $I$ is generated by the orbits of the polynomials in $\mathcal{F}$ under
the action of $\mathfrak{S}_\infty$. 

Given a symmetric ideal $I$, we denote for each $k > 0$, $I_k = \kbb[X_1,\ldots,X_k]$, and 
$V_k(I) = \ZZ(I_k,\kbb^k)$. Clearly, $V_k(I)$  is $\mathfrak{S}_k$-symmetric.

Also, let $\mu = (\mu_1,\ldots,\mu_\ell) \vdash k_0$ be any fixed partition, and for all $k \geq k_0 + \mu_1$, let
\begin{eqnarray}
\label{eqn:def-of-mu-k}
\{\mu\}_k &=& (k-k_0, \mu_1,\mu_2,\ldots,\mu_\ell) \vdash k.
\end{eqnarray} 

It is a consequence of the hook-length formula (Eqn.\eqref{eqn:hook}) that 
\begin{eqnarray}
\label{eqn:polynomial}
\dim_\F(\mathbb{S}^{\{\mu\}_k})=  \frac{\dim_\F(\mathbb{S}_\mu)}{|\mu|!} P_\mu(k),
\end{eqnarray}
where $P_\mu(T)$ is a  monic polynomial having distinct integer roots,  and $\deg(P_\mu) = |\mu|$ (see \cite[7.2.2]{Deligne2004}).

Finally, for a fixed number $p \geq 0$ we pose the following question.

\begin{question}
\label{question:stability}
Let $I \subset A(\R)$ be a finitely generated symmetric ideal.
Does there exist a polynomial $P_{I,p,\mu}(k)$ such that
for all sufficiently large $k$, 
$m_{p,\{\mu\}_k}(V_k(I),\F) = P_{I,p,\mu}(k)$ ? 
In conjunction with
\eqref{eqn:polynomial}, a positive answer would imply that
\[
\dim_\F(\HH^p(V_k(I),\F))_{\{\mu\}_k} =   \frac{\dim_\F(\mathbb{S}_\mu)}{|\mu|!} P_{I,p,\mu}(k) P_\mu(k)
\]
is also given by a polynomial for all large enough $k$.

In particular, taking $\mu=()$ to be the empty partition, is it true that
\[
m_{p,\{\mu\}_k}(V_k(I),\F) = b^p_{\mathfrak{S}_k}(V_k(I),\F)
\]
(that is the $p$-th equivariant Betti number of $V_k(I)$ cf. Notation \ref{not:equivariant-betti} ) is given by a polynomial in $k$ ?

A stronger question is to ask for a bound on the degree of $P_{I,p,\mu}(k)$ as a function of $d,\mu$ and $p$, where $d$ is the maximum of the degrees of the generators of $I$. 
\end{question}

\begin{remark}
Note that it follows from the results of this paper (Theorem \ref{thm:main-product-of-symmetric-quantitative}) that there exists a polynomial  
 $P_{I,p,\mu}(k)$ of degree $O(d^2)$ (where,  $d$ is the maximum of the degrees of the generators of $I$)
with the property that 
 \[
 m_{p,\{\mu\}_k}(V_k(I),\F) \leq  P_{I,p,\mu}(k)
 \]
 for all $k \geq 0$. 
 \end{remark}
 
\begin{remark}
Question \ref{question:stability} has a positive answer for the ideal  $I  \subset A(\R)$, 
generated by the polynomial
\[
f = X_1(X_1-1).
\] 

It is clear from the definition that in this case for each $k > 0$, 
$I_k = (X_1(X_1-1),\ldots,X_k(X_k-1))$, and $V_k(I) = \{0,1\}^k$.

From the discussion in Example \ref{eg:basic} we deduce that for each $p>0,\mu \vdash k_0$, and for all large enough $k$,
\[
m_{p,\{\mu\}_k}(V_k(I),\F) = 0.
\]
For $p=0$, and partitions $\mu$ with $\length(\mu) >1$, we again have for all large enough $k$,
\[
m_{p,\{\mu\}_k}(V_k(I),\F) = 0.
\] 
Finally, for  $p=0$, and any partition $(k_0)$ of length $\leq 1$, and  for all $k \geq  2 k_0$,
\begin{eqnarray*}
m_{0,\{\mu\}_k}(V_k(I),\F) &=& 2(k-k_0) -k+1  \mbox{ (using \eqref{eqn:def-of-mu-k} and \eqref{eqn:even-and-odd})}\\
				      &=& k  - 2k_0 +1.
\end{eqnarray*}
Thus, $m_{p,\{\mu\}_k}(V_k(I),\F)$ is given by a polynomial for all large $k$, for any fixed $p$ and $\mu$.
Notice also that the degree of this polynomial is bounded by $1$.
\end{remark}

\begin{remark}
We point out one crucial difference between
the stability asked for in Question \ref{question:stability} and what is usually meant
by \emph{representational stability} in the FI-module context. In the case of finitely generated FI-modules \cite{Church-et-al}, the multiplicities  $m_{p,\{\mu\}_k}$ are ultimately constant, and in topological applications of the theory, this leads to dimensions of homology groups of each fixed dimension  
(for example, those of the configuration spaces of some fixed manifold) stabilizing to some polynomial (this phenomenon is usually called homological stability). In our case however the multiplicities,
$m_{p,\{\mu\}_k}$,  can grow (albeit polynomially), and the dimensions of the homology groups can grow exponentially 
as seen in Example \ref{eg:basic}.
\end{remark}

\subsection{Algorithmic Conjecture}
As mentioned in the Introduction, a polynomial bound on any topological invariant
of a class of semi-algebraic sets usually implies also that there exists an algorithm with polynomially bounded
complexity for computing it. 
Since we have we proved that the multiplicities of the irreducible representations of 
$\mathfrak{S}_k$ appearing in the cohomology group  of a symmetric $\mathcal{P}$-semi-algebraic set $S \subset \R^k$, where $\deg(P), P \in \mathcal{P}$ is bounded by a constant, is bounded 
by a polynomial function of $\card(\mathcal{P})$ and $k$, the mentioned  principle implies that these multiplicities 
should be computable by an algorithm with polynomially bounded complexity (for fixed $d$).
If this holds, then since the number of irreducibles
that are allowed to appear with positive multiplicity is also polynomially bounded, and
their respective dimensions are polynomially computable using the hook length formula (cf. Eqn.
\eqref{thm:hook}),
we deduce that once these multiplicities are computed, the dimensions of the cohomology groups of $S$ (with coefficients in $\Q$)
can be computed with polynomially bounded complexity. 
 
This leads us to make the following algorithmic conjecture.

\begin{conjecture}
\label{conj:poly}
For any fixed $d >0$, there is an algorithm that takes as input
the description of a symmetric semi-algebraic set $S \subset \R^k$,
defined by a $\mathcal{P}$-closed formula, where 
$\mathcal{P} \subset \R[X_1,\ldots,X_k]^{\mathfrak{S}_k}_{\leq d}$ is a finite set of polynomials, and computes
$m_{i,\lambda}(S,\Q)$, for each $\lambda \vdash k$, and  
$m_{i,\lambda}(S,\Q)>0$, 
as well as all the Betti numbers $b_i(S,\Q)$, with complexity which is polynomial
in  $\card(\mathcal{P})$ and $k$.
\end{conjecture}

\begin{remark}
\label{rem:poly}
We note that Conjecture \ref{conj:poly} is not completely unreasonable, since
an analogous result for computing the generalized Euler-Poincar\'e characteristic
of symmetric semi-algebraic sets has been proved in \cite{BC-conm}. However, computing the
Betti numbers of a semi-algebraic set is usually a much harder task than computing the
Euler-Poincar\'e characteristic. 
More recently, an algorithm with polynomially bounded  complexity has also been given for computing the multiplicities of the trivial representation (i.e. the numbers $m_{i,(k)}(S,\Q)$ using the notation in Conjecture
\ref{conj:poly}) \cite{BC-duke}.
\end{remark}

\bibliographystyle{amsplain}
\bibliography{master}

\def\cprime{$'$} \def\cprime{$'$}
\providecommand{\bysame}{\leavevmode\hbox to3em{\hrulefill}\thinspace}
\providecommand{\MR}{\relax\ifhmode\unskip\space\fi MR }
\providecommand{\MRhref}[2]{%
  \href{http://www.ams.org/mathscinet-getitem?mr=#1}{#2}
}
\providecommand{\href}[2]{#2}
\begin{thebibliography}{10}

\bibitem{Bar97}
A.~I. Barvinok, \emph{On the {B}etti numbers of semialgebraic sets defined by
  few quadratic inequalities}, Math. Z. \textbf{225} (1997), no.~2, 231--244.
  \MR{98f:14044}

\bibitem{Basu1}
S.~Basu, \emph{On bounding the {B}etti numbers and computing the {E}uler
  characteristic of semi-algebraic sets}, Discrete Comput. Geom. \textbf{22}
  (1999), no.~1, 1--18. \MR{1692627 (2000d:14061)}

\bibitem{Bas05-first}
S.~Basu, \emph{Computing the first few {B}etti numbers of semi-algebraic sets
  in single exponential time}, J. Symbolic Comput. \textbf{41} (2006), no.~10,
  1125--1154. \MR{2262087 (2007k:14120)}

\bibitem{Bas05-top}
\bysame, \emph{Computing the top few {B}etti numbers of semi-algebraic sets
  defined by quadratic inequalities in polynomial time}, Found. Comput. Math.
  \textbf{8} (2008), no.~1, 45--80.

\bibitem{BPR8}
S.~Basu, R.~Pollack, and Marie-Fran{\c{c}}oise M.-F.~Roy, \emph{On the {B}etti
  numbers of sign conditions}, Proc. Amer. Math. Soc. \textbf{133} (2005),
  no.~4, 965--974 (electronic). \MR{2117195 (2006a:14096)}

\bibitem{BPR10}
S.~Basu, R.~Pollack, and M.-F. Roy, \emph{Betti number bounds, applications and
  algorithms}, Current Trends in Combinatorial and Computational Geometry:
  Papers from the Special Program at MSRI, MSRI Publications, vol.~52,
  Cambridge University Press, 2005, pp.~87--97.

\bibitem{BPRbook2}
\bysame, \emph{Algorithms in real algebraic geometry}, Algorithms and
  Computation in Mathematics, vol.~10, Springer-Verlag, Berlin, 2006 (second
  edition). \MR{1998147 (2004g:14064)}

\bibitem{BPRbettione}
\bysame, \emph{Computing the first {B}etti number of a semi-algebraic set},
  Found. Comput. Math. \textbf{8} (2008), no.~1, 97--136.

\bibitem{BC2013}
S.~{Basu} and C.~{Riener}, \emph{{Bounding the equivariant Betti numbers and
  computing the generalized Euler-Poincar\'{e} characteristic of symmetric
  semi-algebraic sets}}, ArXiv e-prints (2013).

\bibitem{BC-conm}
\bysame, \emph{{Efficient algorithms for computing the Euler-Poincar\'e
  characteristic of symmetric semi-algebraic sets}}, ArXiv e-prints (to appear
  in AMS Contemporary Mathematics) (2016).

\bibitem{BC-duke}
\bysame, \emph{{On the equivariant Betti numbers of symmetric semi-algebraic
  sets: vanishing, bounds and algorithms}}, ArXiv e-prints (2016).

\bibitem{Basu-sheaf}
Saugata Basu, \emph{A {C}omplexity {T}heory of {C}onstructible {F}unctions and
  {S}heaves}, Found. Comput. Math. \textbf{15} (2015), no.~1, 199--279.
  \MR{3303696}

\bibitem{BP'R07joa}
Saugata Basu, Dmitrii~V. Pasechnik, and Marie-Fran{\c{c}}oise Roy,
  \emph{Computing the {B}etti numbers of semi-algebraic sets defined by partly
  quadratic sytems of polynomials}, J. Algebra \textbf{321} (2009), no.~8,
  2206--2229. \MR{2501518 (2010a:14092)}

\bibitem{BC-advances}
Saugata Basu and Cordian Riener, \emph{Bounding the equivariant {B}etti numbers
  of symmetric semi-algebraic sets}, Adv. Math. \textbf{305} (2017), 803--855.
  \MR{3570148}

\bibitem{Ceccherini-book}
Tullio Ceccherini-Silberstein, Fabio Scarabotti, and Filippo Tolli,
  \emph{Representation theory of the symmetric groups}, Cambridge Studies in
  Advanced Mathematics, vol. 121, Cambridge University Press, Cambridge, 2010,
  The Okounkov-Vershik approach, character formulas, and partition algebras.
  \MR{2643487 (2011h:20024)}

\bibitem{Church-et-al}
Thomas Church, Jordan~S. Ellenberg, and Benson Farb, \emph{F{I}-modules and
  stability for representations of symmetric groups}, Duke Math. J.
  \textbf{164} (2015), no.~9, 1833--1910. \MR{3357185}

\bibitem{Dalbec}
John Dalbec, \emph{Multisymmetric functions}, Beitr\"age Algebra Geom.
  \textbf{40} (1999), no.~1, 27--51. \MR{1678567 (2000f:05084)}

\bibitem{Deligne2004}
P.~Deligne, \emph{La cat\'egorie des repr\'esentations du groupe sym\'etrique
  {$S_t$}, lorsque {$t$} n'est pas un entier naturel}, Algebraic groups and
  homogeneous spaces, Tata Inst. Fund. Res. Stud. Math., Tata Inst. Fund. Res.,
  Mumbai, 2007, pp.~209--273. \MR{2348906 (2009b:20021)}

\bibitem{GVZ04}
A.~Gabrielov, N.~Vorobjov, and T.~Zell, \emph{Betti numbers of semialgebraic
  and sub-{P}faffian sets}, J. London Math. Soc. (2) \textbf{69} (2004), no.~1,
  27--43. \MR{2025325 (2004k:14105)}

\bibitem{GV07}
Andrei Gabrielov and Nicolai Vorobjov, \emph{Approximation of definable sets by
  compact families, and upper bounds on homotopy and homology}, J. Lond. Math.
  Soc. (2) \textbf{80} (2009), no.~1, 35--54. \MR{2520376}

\bibitem{GRW2016}
Paul G\"orlach, Cordian Riener, and Tillmann Wei\ss~er, \emph{Deciding
  positivity of multisymmetric polynomials}, J. Symbolic Comput. \textbf{74}
  (2016), 603--616. \MR{3424059}

\bibitem{Hirzebruch-book}
Friedrich Hirzebruch, \emph{Topological methods in algebraic geometry},
  Classics in Mathematics, Springer-Verlag, Berlin, 1995, Translated from the
  German and Appendix One by R. L. E. Schwarzenberger, With a preface to the
  third English edition by the author and Schwarzenberger, Appendix Two by A.
  Borel, Reprint of the 1978 edition. \MR{1335917 (96c:57002)}

\bibitem{Milnor2}
J.~Milnor, \emph{On the {B}etti numbers of real varieties}, Proc. Amer. Math.
  Soc. \textbf{15} (1964), 275--280. \MR{0161339 (28 \#4547)}

\bibitem{Mimura-Toda-book}
Mamoru Mimura and Hirosi Toda, \emph{Topology of {L}ie groups. {I}, {II}},
  Translations of Mathematical Monographs, vol.~91, American Mathematical
  Society, Providence, RI, 1991, Translated from the 1978 Japanese edition by
  the authors. \MR{1122592 (92h:55001)}

\bibitem{OP}
I.~G. Petrovski{\u\i} and O.~A. Ole{\u\i}nik, \emph{On the topology of real
  algebraic surfaces}, Izvestiya Akad. Nauk SSSR. Ser. Mat. \textbf{13} (1949),
  389--402. \MR{0034600 (11,613h)}

\bibitem{Procesi-book}
Claudio Procesi, \emph{Lie groups}, Universitext, Springer, New York, 2007, An
  approach through invariants and representations. \MR{2265844 (2007j:22016)}

\bibitem{Riener}
Cordian Riener, \emph{On the degree and half-degree principle for symmetric
  polynomials}, J. Pure Appl. Algebra \textbf{216} (2012), no.~4, 850--856.
  \MR{2864859}

\bibitem{Schapira-notes}
Pierre Schapira, \emph{Algebra and topology}, Course at Paris VI University,
  2007/2008.

\bibitem{Scheiblechner07}
Peter Scheiblechner, \emph{On the complexity of deciding connectedness and
  computing {B}etti numbers of a complex algebraic variety}, J. Complexity
  \textbf{23} (2007), no.~3, 359--379. \MR{2330991 (2009d:14020)}

\bibitem{Stanley-vol1}
Richard~P. Stanley, \emph{Enumerative combinatorics. {V}olume 1}, second ed.,
  Cambridge Studies in Advanced Mathematics, vol.~49, Cambridge University
  Press, Cambridge, 2012. \MR{2868112}

\bibitem{T}
R.~Thom, \emph{Sur l'homologie des vari\'et\'es alg\'ebriques r\'eelles},
  Differential and Combinatorial Topology (A Symposium in Honor of Marston
  Morse), Princeton Univ. Press, Princeton, N.J., 1965, pp.~255--265.
  \MR{0200942 (34 \#828)}

\bibitem{Timofte03}
Vlad Timofte, \emph{On the positivity of symmetric polynomial functions. {I}.
  {G}eneral results}, J. Math. Anal. Appl. \textbf{284} (2003), no.~1,
  174--190. \MR{1996126 (2005c:05196)}

\bibitem{Walther1}
Uli Walther, \emph{Algorithmic determination of the rational cohomology of
  complex varieties via differential forms}, Symbolic computation: solving
  equations in algebra, geometry, and engineering (South Hadley, MA, 2000),
  Contemp. Math., vol. 286, Amer. Math. Soc., Providence, RI, 2001,
  pp.~185--206. \MR{1874280 (2003b:14027)}

\end{thebibliography}

\end{document}